\documentclass{amsart}

\usepackage{amsmath,amssymb,amsfonts,amssymb,amsthm}
\usepackage{verbatim}
\usepackage[usenames]{color}
\usepackage{hyperref}
\usepackage{url}

\usepackage{graphicx}

\usepackage{amsthm,graphicx,appendix}

\usepackage{tikz,tikz-qtree,ifthen,cancel}
\usepackage{array}
\usepackage{adjustbox}
\usepackage{appendix}
\usetikzlibrary{calc,shapes,patterns,positioning}

\numberwithin{equation}{section}

\theoremstyle{plain}
\newtheorem{theorem}{Theorem}[section]

\theoremstyle{definition}
\newtheorem{definition}[theorem]{Definition}

\newtheorem{question}[theorem]{Question}
\newtheorem{problem}[theorem]{Problem}

\newcommand{\Erdos}{Erd{\H{o}}s}
\newcommand{\Fraisse}{Fra{\"{i}}ss{\'{e}}}

\newcommand{\Hubicka}{Hubi{\v{c}}ka}
\newcommand{\Lauchli}{L{\"{a}}uchli}
\newcommand{\Masulovic}{Ma{\v{s}}ulovi{\'{c}}}
\newcommand{\Nesetril}{Ne{\v{s}}et{\v{r}}il}

\newcommand{\Rodl}{R{\"{o}}dl}

\newcommand{\om}{\omega}

\newcommand{\bfA}{\mathbf{A}}
\newcommand{\bfB}{\mathbf{B}}
\newcommand{\bfC}{\mathbf{C}}
\newcommand{\bfD}{\mathbf{D}}
\newcommand{\bfE}{\mathbf{E}}

\newcommand{\bK}{\mathbf{K}}
\newcommand{\bR}{\mathbf{R}}
\newcommand{\bfS}{\mathbf{S}}

\newcommand{\bQ}{\mathbb{Q}}
\newcommand{\bS}{\mathbb{S}}

\newcommand{\sse}{\subseteq}
\newcommand{\re}{\upharpoonright}

\newcommand{\ra}{\rightarrow}

\newcommand{\trl}{\triangleleft}

\newcommand{\lgl}{\langle}
\newcommand{\rgl}{\rangle}

\DeclareMathOperator{\dom}{dom}
\DeclareMathOperator{\height}{ht}\DeclareMathOperator{\Succ}{Succ}
\DeclareMathOperator{\mc}{cl}
\DeclareMathOperator{\type}{tp}
\DeclareMathOperator{\Age}{Age}

\begin{document}

%------
% Insert the title of your paper and (if necessary)
% a short title for the running head.
%------
\title[Ramsey theory of homogeneous structures]{Ramsey theory of homogeneous structures: current trends and open problems}

\author{Natasha Dobrinen}
\address{Department of Mathematics\\
 University of Denver \\
2390 S.\ York St.\\ Denver, CO \ 80208 U.S.A.}
\email{natasha.dobrinen@du.edu}
  \urladdr{\url{https://cs.du.edu/dobrinen/}}
\thanks{The author gratefully acknowledges  support from   National Science Foundation Grant DMS-1901753}

%------
% Insert your abstract.
%------
\begin{abstract}

%The abstract should not be more than 10 lines long.

This article
highlights  historical achievements in the  partition theory of countable homogeneous relational structures, and presents recent work, current trends, and open problems.
Exciting   recent developments include  new methods involving logic, topological Ramsey spaces, and  category theory.
The paper concentrates on
 big Ramsey degrees, presenting their essential structure where known and outlining areas for further development.
Cognate areas, including infinite dimensional Ramsey theory of homogeneous structures and  partition theory of uncountable structures, are also discussed.
\end{abstract}

\maketitle

 This article will appear in the 2022 ICM Proceedings.
It  is dedicated to Norbert Sauer for his seminal works on the  partition theory of homogeneous structures, and  for his mathematical and personal generosity.

%%%%%%%%%%%%%%%%%%%%%%%%%%%%%%%%
%%%%%%%%%%%%%%%%%%%%%%%%%%%%%%%%
%%%%%%%%%%%%%%%%%%%%%%%%%%%%%%%%

\section{Introduction}

%REMEMBER WHO MY AUDIENCE IS - THE GENERAL MATHEMATICIAN WHOM I WANT TO CONVINCE THAT THIS IS A TOTALLY COOL SUBJECT.

%DON'T TRY TO BE TOO FANCY.  GET WHAT I WANT TO SAY OUT TO THE AUDIENCE IN AS CLEAR AND SIMPLE A MANNER AS POSSIBLE.

%ONE PARAGRAPH OF INTRO AND CONNECTIONS

%THIS IS AN IDEAS ARTICLE.  I DON'T HAVE TO GIVE THE DETAILS.  I DO NEED TO GIVE THE BIG PICTURE AND MY INSIGHTS AS TO WHAT IS REALLY GOING ON.

Ramsey theory  is a beautiful
 subject  which interrelates with   a multitude of  mathematical
fields.
In particular,
since its inception, developments in Ramsey theory have often been motivated by problems in logic;
 in turn, Ramsey theory has instigated some seminal developments in logic.
The intent of this article is to provide the general mathematician with an introduction to the
 intriguing subject of Ramsey theory on homogeneous  structures while
being detailed enough to
describe the state of the art and
 the main ideas at play.
We
 present
  historical highlights and  discuss
why
solutions to problems on homogeneous structures
require more than just straightforward applications of
 finite structural  Ramsey theory.
In the following sections, we
map out  collections of recent results and
   methods which were
 developed to overcome  obstacles associated with forbidden substructures.
These new methods involve applications from  logic (especially forcing but also ideas from model theory),  topological Ramsey spaces, and category theory.
%The article will close with future directions and   open problems.

The subject of Ramsey theory on infinite structures begins with this lovely theorem.

\begin{theorem}[Ramsey, \cite{Ramsey30}]\label{thm.RamseyInfinite}
Given positive integers $k$  and $r$ and a coloring
of the $k$-element subsets of the natural numbers
$\mathbb{N}$ into $r$ colors,
there is an infinite set of natural numbers $N\sse\mathbb{N}$ such that all $k$-element subsets of $N$ have the same color.
\end{theorem}

There are two natural interpretations of
Ramsey's Theorem in terms of infinite  structures.
First,  letting $<$ denote the standard linear order on $\mathbb{N}$,
 Ramsey's Theorem shows that given any finite  coloring of all
linearly ordered substructures of $(\mathbb{N},<)$ of size $k$,
there is an  isomorphic
substructure $(N,<)$ of $(\mathbb{N},<)$
such that  all  linearly ordered substructures of  $(N,<)$  of size $k$ have the same color.
Second, one may think of the $k$-element subsets of
$\mathbb{N}$ as $k$-hyperedges.
Then Ramsey's Theorem yields that, given any finite coloring of the $k$-hyperedges  of  the complete $k$-regular hypergraph on infinitely many vertices, there is
an isomorphic subgraph
in which all $k$-hyperedges have the same color.

Given this, one might naturally wonder about other structures.

\begin{question}\label{q.beginners}
Which infinite structures carry an analogue of Ramsey's Theorem?
\end{question}

The rational numbers
 $(\mathbb{Q},<)$ as   a dense linearly ordered structure (without endpoints)
was  the earliest test case.
It is a fun exercise to show that  given any coloring of the rational numbers into finitely many colors, there is one color-class  which contains a dense linear order, that is,  an isomorphic subcopy of the rationals in one color.
Thus,   the rationals satisfy
a structural  pigeonhole principle known as   {\em indivisibility}.

The direct analogy with Ramsey's theorem ends, however,  when we consider pairs of rationals.
It follows from work of Sierpi\'{n}ski
 in \cite{Sierpinski}
 that there is a coloring of the pairs of rationals into  two colors so that
both colors persist
in every isomorphic subcopy of the rationals.
Sierpi\'{n}ski's coloring  provides a clear
 understanding of one of the fundamental issues arising
in partition theory of infinite structures
 not occurring in
 finite structural Ramsey theory.
Let   $\{q_i:i\in \mathbb{N}\}$ be a listing of  the rational numbers, without repetition, and
for $i<j$ define
$c(\{q_i,q_j\}) =
          \mathrm{blue}$ if  $ q_i<q_j$,
 and
         $ c(\{q_i,q_j\}) =  \mathrm{red}$  if   $q_j<q_i$.
Then in each subset $Q\subseteq \mathbb{Q}$ forming  a dense linear order,
both color classes persist; that is,
there are pairs  of rationals in $Q$ colored red and also pairs of rationals in $Q$ colored blue.
Since it is impossible to find an isomorphic subcopy of the rationals in which all pairsets have the same color,
a direct analogue of Ramsey's theorem does not hold for the rationals.

The failure of  the straightforward  analogue of Ramsey's theorem  is not the end, but rather  just the  beginning of the story.
  Galvin (unpublished)
showed a few decades later that
there is  a bound on the number of unavoidable colors:
Given any  coloring of the pairs of rationals into finitely many colors, there is a subcopy  of the rationals
 in which all pairs  belong to the union of two color classes.
Now one sees that Question \ref{q.beginners} ought to  be refined.

\begin{question}\label{q.refined}
For which infinite structures $\bfS$ is there a Ramsey-analogue in the following sense:
Let $\bfA$ be a finite substructure of $\bfS$.
Is there a positive integer $T$  such that for  any coloring
of the copies of $\bfA$ into finitely many colors,
there is a subcopy $\bfS'$ of $\bfS$
in which there are no more than
$T$ many colors for the copies of $\bfA$?
\end{question}

The least such integer $T$, when it exists,
is denoted $T(\bfA)$ and
called the {\em big Ramsey degree} of $\bfA$ in $\bfS$,
 a term coined in Kechris--Pestov--Todorcevic (2005).
The ``big'' refers to the fact that
we require  an isomorphic subcopy of an
{\em infinite} structure in which the number of colors is as small as possible
(in contrast to the concept of {\em small Ramsey degree} in finite structural Ramsey theory).

Notice how Sierpi\'{n}ski played
the  enumeration $\{q_i:i\in\mathbb{N}\}$
 of the rationals  against
the dense linear order to construct a coloring of pairsets of rationals into two colors, each of which persists in every subcopy of the rationals.
This simple but deep idea
sheds light on a  fundamental
difference between
 finite  and infinite structural Ramsey theory.
The interplay between the enumeration and the
relations on an infinite structure
has  bearing on the number of colors that must persist in any subcopy of that  structure.
We will see examples of this at work throughout this article and explain the  general principles which have been found for certain classes of  structures with relations of arity at most two, even  as the subject aims towards a future  overarching theory of big Ramsey degrees.

%%%%%%%%%%%%%%%%
%%%%%%%%%%%%%%%%
%%%%%%%%%%%%%%%%
%%%%%%%%%%%%%%%%

\section{The Questions}\label{sec.theproblems}

%We  now set forth  the driving questions in the Ramsey theory of homogenous structures. Minimal background needed for   a clear understanding of these questions is provided.

Given a finite
relational
language
 $\mathcal{L}=\{R_i: i< k\}$ with each
relation symbol  $R_i$ of some finite arity, say $n_i$,
an {\em $\mathcal{L}$-structure}
is a tuple
 $\bfA=\lgl A, R_0^{\bfA},\dots, R_{k-1}^{\bfA}\rgl$,
where $A \ne\emptyset$ is the {\em universe} of $\bfA$  and for each $i<k$, $R_i^{\bfA}\sse A^{n_i}$.
For $\mathcal{L}$-structures $\bfA$ and $\bfB$,
an {\em embedding} from
$\bfA$ into $\bfB$  is an injection $e:A\ra B$ such that
for all $i<k $, $R_i^{\bfA}(a_1,\dots,a_{n_i})\leftrightarrow R_i^{\bfB}(e(a_1),\dots,e(a_{n_i}))$.
The $e$-image of $\bfA$ is  a {\em copy} of $\bfA$ in $\bfB$.
If $e$ is the identity map, then
 $\bfA$ is a {\em substructure} of $\bfB$.
An  {\em isomorphism}  is an  embedding which is onto its image.
We write $\bfA\le\bfB$  exactly when  there is an embedding of  $\bfA$  into $\bfB$, and
 $\bfA\cong\bfB$  exactly when  $\bfA$ and $\bfB$ are isomorphic.

A class $\mathcal{K}$ of finite structures for a relational language $\mathcal{L}$
 is called a {\em \Fraisse\ class}  if it is hereditary, satisfies the joint embedding and amalgamation properties, contains (up to isomorphism) only countably many structures, and contains structures of arbitrarily large finite cardinality.
 $\mathcal{K}$ is  {\em hereditary} if whenever $\bfB\in\mathcal{K}$ and  $\bfA\le\bfB$, then also $\bfA\in\mathcal{K}$.
$\mathcal{K}$ satisfies the {\em joint embedding property} if for any $\bfA,\bfB\in\mathcal{K}$,
there is a $\bfC\in\mathcal{K}$ such that $\bfA\le\bfC$ and $\bfB\le\bfC$.
 $\mathcal{K}$ satisfies the {\em amalgamation property} if for any embeddings
$f:\bfA\ra\bfB$ and $g:\bfA\ra\bfC$, with $\bfA,\bfB,\bfC\in\mathcal{K}$,
there is a $\bfD\in\mathcal{K}$ and  there are embeddings $r:\bfB\ra\bfD$ and $s:\bfC\ra\bfD$ such that
$r\circ f = s\circ g$.
A \Fraisse\ class $\mathcal{K}$ satisfies the {\em strong amalgamation  property} (SAP) if given $\bfA,\bfB,\bfC\in\mathcal{K}$ and  embeddings $e:\bfA\ra \bfB$ and $f:\bfA\ra \bfC$,
there is some $\bfD\in\mathcal{K}$ and embeddings $e':\bfB\ra\bfD$ and $f':\bfC\ra\bfD$ such that $e'\circ e=f'\circ f$, and $e'[B]\cap f'[C]=e'\circ e[A]=f'\circ f[A]$.
We say that  $\mathcal{K}$ satisfies the {\em free amalgamation property} (FAP) if it satisfies the SAP and moreover, $\bfD$ can be chosen so that $\bfD$ has no additional relations other than those inherited from $\bfB$ and $\bfC$.

Let $\bfA,\bfB,\bfC$ be   $\mathcal{L}$-structures
such that $\bfA\le\bfB\le \bfC$.
We use ${\bfB\choose\bfA}$ to denote the set of all
copies of $\bfA$ in
 $\bfB$.
The \Erdos--Rado arrow notation
 $\bfC\ra(\bfB)_k^{\bfA}$
means  that for each coloring of ${\bfC\choose \bfA}$ into $k$ colors, there is a $\bfB' \in {\bfC\choose\bfB}$ such that
${\bfB'\choose\bfA}$ is  {\em monochromatic}, meaning
 every member of ${\bfB'\choose\bfA}$ has the same color.

 \begin{definition}\label{defn.RP}
A \Fraisse\ class  $\mathcal{K}$ has the {\em Ramsey property} if  for any two structures $\bfA\le\bfB$ in $\mathcal{K}$ and any  $k\ge 2$,
there is a $\bfC\in\mathcal{K}$ with $\bfB\le\bfC$ such that
$\bfC\ra (\bfB)^{\bfA}_k$.
\end{definition}

 %Ramsey theory of \Fraisse\ classes is its own glorious subject to which we cannot do justice within the scope of this  paper.
Many  \Fraisse\ classes, such as the class of finite graphs, do not have the Ramsey property.
However,   by allowing a finite expansion of the language,  often by just a  linear order,
the Ramsey property  becomes more feasible.
Letting $<$ be a binary relation symbol not in the language $\mathcal{L}$ of
$\mathcal{K}$,
an $\mathcal{L}\cup\{<\}$-structure  is in
 $\mathcal{K}^{<}$  if and only if its universe is linearly ordered by $<$ and its $\mathcal{L}$-reduct  is a member of $\mathcal{K}$.
A highlight is
  the work of \Nesetril\ and \Rodl\ in
\cite{Nesetril/Rodl77} and \cite{Nesetril/Rodl83},
proving that for any   \Fraisse\ class $\mathcal{K}$  with FAP, its ordered version $\mathcal{K}^{<}$
 has the Ramsey property.
The  recent paper \cite{HN19} by \Hubicka\ and \Nesetril\
presents the state of the art in finite structural Ramsey theory.
 Examples of \Fraisse\ classes with the Ramsey property include the class of  finite linear orders, and the classes of finite ordered versions of  graphs, digraphs, tournaments, triangle-free graphs, posets,
metric spaces,
hypergraphs, hypergraphs omitting some irreducible substructures, and many more.

A structure $\bK$ is called {\em universal} for a class of structures $\mathcal{K}$ if each member of $\mathcal{K}$ embeds into $\bK$.
A structure $\bK$ is {\em homogeneous} if  each isomorphism between finite substructures of $\bK$ extends to an automorphism of $\bK$.
Unless otherwise specified, we will
 write
 {\em homogeneous} to  mean {\em countably infinite homogeneous}, such structures being  the focus of this paper.
The {\em age}
 of an infinite structure $\bK$, denoted  Age$(\bK)$, is the collection of all finite structures which embed into $\bK$.
A fundamental theorem of \Fraisse\ from   \cite{Fraisse54} shows that
each \Fraisse\ class gives rise to a homogeneous structure via a construction called the {\em \Fraisse\ limit}.
Conversely,
given any countable  homogeneous structure $\bK$,
Age$(\bK)$ is a \Fraisse\ class and, moreover,
the \Fraisse\ limit of Age$(\bK)$ is isomorphic to $\bK$.
The Kechris--Pestov--Todorcevic correspondence between
the Ramsey property of a \Fraisse\ class
and
extreme amenability of the automorphism group of its \Fraisse\ limit in \cite{KPT05}
 propelled  a
burst of discoveries of more \Fraisse\ classes with the Ramsey property.

First an esoteric but driving question in the area.

\begin{question}\label{q.essence}
What {\em is}   a big Ramsey degree?
\end{question}

What is the essential nature of a big Ramsey degree?
Why is it that given a \Fraisse\ class $\mathcal{K}$ satisfying the Ramsey property,
its \Fraisse\ limit  usually fails to carry the full analogue of Ramsey's Theorem \ref{thm.RamseyInfinite} (i.e.\ all big Ramsey degrees being one)?
A theorem of Hjorth in \cite{Hjorth}
showed that
for any homogeneous structure $\bK$ with
$|\mathrm{Aut}(\bK)|>1$,
there is a structure in Age$(\bK)$ with big Ramsey degree at least two.
While much remains open, we now
 have an answer to Question \ref{q.essence}   for  FAP  and some SAP homogeneous structures with   finitely many  relations of arity at most two, and these
 results  will be discussed  in  the following sections.

We say that $\bfS$ has {\em finite big Ramsey degrees} if $T(\bfA)$ exists for each finite substructure $\bfA$ of $\bfS$.
We say that {\em exact big Ramsey degrees} are known if there is either a computation of the degrees or a characterization from which they can be computed.
 {\em Indivisibility} holds if
$T(\bfA)=1$  for each  one-element  substructure  $\bfA$ of $\bfS$.
The following questions  progress in order of strength: A positive answer to (3) implies a positive answer to (2), which in turn implies a positive answer to (1).

\begin{question}\label{q.3}
Given a homogeneous structure $\bK$,
\begin{enumerate}
\item
Does $\bK$ have finite big Ramsey degrees?
That is, can one find upper bounds ensuring that  big Ramsey degrees exist?
\item
If $\bK$ has finite big Ramsey degrees,
is there a
 characterization of  the exact
 big Ramsey degrees via canonical partitions?
If yes,  calculate or find an algorithm to calculate them.
\item
Does $\bK$ carry a big Ramsey structure?
\end{enumerate}
\end{question}

Part (2) of this question  involves finding  {\em canonical partitions}.

\begin{definition}[Canonical Partition, \cite{LSV06}]\label{defn.cp}
Given a \Fraisse\ class
$\mathcal{K}$  with \Fraisse\ limit $\bK$, and given
$\bfA\in\mathcal{K}$,
a partition
$\{P_i:i<n\}$
 of ${\bK\choose \bfA}$ is  {\em canonical}
if the following hold:
For each finite coloring  of ${\bK\choose \bfA}$, there is a subcopy
$\bK'$ of $\bK$
such that
for each $i<n$, all members of $P_i\cap {\bK'\choose\bfA}$
have the  same  color; and
{\em persistence}:
For every subcopy
$\bK'$ of $\bK$
and each $i < n$,
 $P_i\cap {\bK'\choose\bfA}$ is non-empty.
\end{definition}

Canonical partitions recover an exact analogue of Ramsey's theorem for each piece of the partition.
In practice such partitions are characterized by adding extra structure to $\bK$,
 including the enumeration of the universe of $\bK$ and a tree-like structure capturing the relations of $\bK$ against the enumeration.

Part (3) of Question \ref{q.3} has to do with
a connection between big Ramsey degrees and topological dynamics, in the spirit of the  Kechris--Pestov--Todorcevic correspondence,
proved by
 Zucker in  \cite{Zucker19}.
A {\em big Ramsey structure} is essentially  a finite expansion $\bK^*$ of $\bK$ so that each finite substructure of $\bK^*$ has big Ramsey degree one, and moreover the unavoidable colorings cohere in that for $\bfA,\bfB\in\mathrm{Age}(\bK)$ with $\bfA$ embedding into $\bfB$,
 the canonical partition for copies of $\bfB$  when restricted to copies of $\bfA$ recovers the canonical partition for copies of $\bfA$.
Big Ramsey structures imply canonical partitions.
The reverse is not known in general, but
  certain types of canonical partitions are known to  imply big Ramsey structures  (Theorem 6.10 in  \cite{CDP21}), and it seems reasonable to the author to expect that  (1)--(3) are equivalent.

Canonical partitions and big Ramsey structures are really getting at the question of whether we can  find an optimal finite expansion $\bK^*$ of  a given homogeneous  structure  $\bK$
so that $\bK^*$  carries an exact analogue of Ramsey's Theorem.
In this sense,
 big Ramsey degrees are not quite so mysterious, but  are rather saying that an exact analogue of Ramsey's theorem holds
  for an appropriately expanded structure.
The question then becomes, what is the appropriate expansion?

%%%%%%%%%%%%%%%%%%%%%%
%%%%%%%%%%%%%%%%%%%%%%
%%%%%%%%%%%%%%%%%%%%%%
%%%%%%%%%%%%%%%%%%%%%%

\section{Case Study:  The rationals}\label{sec.CaseStudies}

The  big Ramsey degrees for the rationals
were determined  by 1979.
Laver in 1969 (unpublished, see \cite{DevlinThesis})
 utilized a Ramsey theorem for trees due to Milliken \cite{Milliken79} (Theorem \ref{thm.Milliken}) to find upper bounds.
Devlin  completed the picture in his Ph.D.\ thesis \cite{DevlinThesis}, calculating the big Ramsey degrees of the rationals.
These surprisingly turn out to be related  to  the odd  coefficients in  the Taylor series of the tangent function: The big Ramsey degree for $n$-element subsets of the rationals
is
 $T(n)=(2n-1)!c_{2n-1}$, where
$c_k$ is the  $k$-th coefficient in the Taylor series for the tangent function,
$\tan(x)=\sum_{k=0}^{\infty}c_k x^k$.
As Todorcevic states,
the big Ramsey degrees for the rationals
 ``characterize the Ramsey theoretic properties of the countable dense linear ordering $(\mathbb{Q},<)$ in a very precise sense.
The numbers $T(n)$ are some sort of Ramsey degrees that measure the complexity of an arbitrary finite coloring of the $n$-element subsets of $\mathbb{Q}$  modulo, of course, restricting to
the $n$-element subsets of $X$ for some appropriately chosen dense linear subordering $X$ of $\mathbb{Q}$.''
 (page 143, \cite{TodorcevicBK10}, notation modified)

We    present Devlin's characterization of the big Ramsey degrees of the rationals and the four main steps in his proof.
(A detailed proof  appears in  Section 6.3 of \cite{TodorcevicBK10}.)
Then  we will present a  method  from \cite{CDP21} using coding trees of $1$-types which
bypasses
 non-essential constructs, providing what we see as a satisfactory answer to Question \ref{q.essence} for the rationals.

We use some standard  mathematical logic notation, providing definitions as needed for the general mathematician.
The  set of all natural numbers $\{0,1,2,\dots\}$
is denoted by $\om$.
Each natural number
$k\in\om$  is equated with  the set
$\{0,\dots, k-1\}$ and  its natural linear ordering.
 $k\in\om$ and $k<\om$  are synonymous.
For $k\in \om$,
  $k^{<\om}$ denotes the tree of all finite sequences with entries in $\{0,\dots, k-1\}$,
and
$\om^{<\om}$ denotes the tree of all finite sequences of natural numbers.
Finite sequences with any sort of entries are thought of as functions with domain some natural number.
Thus, for  a finite sequence $t$
 the {\em length} of $t$, denoted  $|t|$, is the domain of the function $t$, and for $i\in\dom(t)$, $t(i)$ denotes the $i$-th entry of the sequence $t$.
For $\ell\in\om$, we write $t\re \ell$ to denote the initial segment of $t$ of length $\ell$ if $\ell\le|t|$, and $t$ otherwise.
For two finite sequences $s$ and $t$,
we write $s\sqsubseteq t$ when  $s$ is an initial segment of $t$,
%meaning that $t\re |s|=s$,
and we write  $s\sqsubset t$   when
$s$ is a proper initial segment of $t$, meaning that
$s\sqsubseteq t$ and $s\ne t$.
We write $s\wedge t$ to denote
the {\em meet} of $s$ and $t$; that is,
the longest sequence  which is an initial segment of both $s$ and $t$.
Given a subset $S$ of a tree of finite sequences,
the {\em meet closure} of $S$, denoted $\mc(S)$,
is the set of all nodes in $S$ along with the set of all meets $s\wedge t$, for $s,t\in S$.

A Ramsey theorem for trees, due to Milliken, played a central role in Devlin's work and  has informed   subsequent approaches to finding upper bounds for big Ramsey degrees.
In this area,
 a subset $T\sse\om^{<\om}$ is called  a {\em tree} if there is a subset $L_T\sse\om$  such that
$T=\{t\re \ell:t\in T,\ \ell\in L_T\}$.
Thus,  a tree is closed under initial segments of lengths in $L_T$,
 but not necessarily closed under all initial segments in $\om^{<\om}$.
The  {\em  height} of a node  $t$ in $T$,  denoted $\height_T(t)$, is the order-type of the set  $\{s\in T:s\sqsubset  t\}$, linearly ordered by $\sqsubset$.
We write
 $T(n)$ to  denote $\{t\in T:\height_T(t)=n\}$.
 For $t\in T$,  let $\Succ_T(t)=\{s\re(|t|+1):s\in T$ and $t\sqsubset s\}$, noting that $\Succ_T(t)\sse T$     only if $|t|+1\in L_T$.

A subtree $S\sse T$
 is a {\em strong subtree} of $T$ if
$L_S\sse L_T$
 and each node $s$ in $S$ branches
 as widely as $T$ will allow,
meaning that for $s\in S$,
for each
 $t \in \Succ_T(s)$
there is an extension $s'\in S$ such that
$t\sqsubseteq s'$.
For the next theorem, define
$ \prod_{i<d} T_i(n)$ to be the set of sequences  $(t_0,\dots,t_{d-1})$ where $t_i\in T_i(n)$, the product of the $n$-th levels of the trees $T_i$.
 Then let
\begin{equation}
\bigotimes_{i<d}T_i:=\bigcup_{n<\om}  \prod_{i<d} T_i(n).
\end{equation}

The following is the strong tree version of the Halpern--\Lauchli\ Theorem.

\begin{theorem}[Halpern--\Lauchli, \cite{HL66}]\label{thm.HL}
Let   $d$ be a positive integer,
 $T_i\sse \om^{<\om}$  ($i<d$) be  finitely branching trees
 with no terminal nodes, and  $r\ge 2$.
Given a coloring
$c:\bigotimes_{i<d} T_i \ra r$,
there is an increasing sequence $\lgl m_n:n<\om\rgl$ and
 strong  subtrees $S_i\le T_i$
 such that for all $i<d$ and $n<\om$, $S_i(n)\sse T_i(m_n)$, and
$c$ is constant on $\bigotimes_{i<d}S_i$.
\end{theorem}

The Halpern--\Lauchli\ Theorem  has  a particularly strong  connection with logic.
It  was isolated by Halpern and  L\'{e}vy as a key juncture in their work to prove that
the Boolean Prime Ideal Theorem is strictly weaker than the Axiom of Choice  over
the Zermelo--Fraenkel Axioms  of set theory.
Once proved by Halpern and \Lauchli,
Halpern and L\'{e}vy completed their proof in  \cite{Halpern/Levy71}.

Harrington (unpublished)
devised an innovative  proof of the
Halpern--\Lauchli\ Theorem which used Cohen forcing.
The forcing helps find good nodes in the trees $T_i$ from which to start building
the subtrees $S_i$.
From then on, the forcing is used $\om$ many times, each time
 running  an unbounded search for  finite  sets $S_i(n)$   which  satisfy that level of the Halpern--\Lauchli\ Theorem.
Being finite, each $S_i(n)$ is  in the ground model.
The proof   entails  neither passing to a generic extension nor any use of  Shoenfield's Absoluteness Theorem.

A {\em $k$-strong subtree} is a strong subtree with $k$ many levels.
The following theorem  is
proved  inductively using  Theorem \ref{thm.HL}.

\begin{theorem}[Milliken, \cite{Milliken79}]\label{thm.Milliken}
Let $T\sse \om^{<\om}$ be a finitely branching tree with no terminal nodes,
  $k\ge 1$, and  $r\ge 2$.
Given a coloring
of all
$k$-strong subtrees of
 $T$ into $r$ colors,
 there is an infinite  strong subtree $S\sse T$ such that all  $k$-strong subtrees of  $S$ have the same color.
\end{theorem}

For more on the Halpern--\Lauchli\ and Milliken Theorems, see \cite{DKBK},  \cite{Larson12}, and  \cite{TodorcevicBK10}.
Now we look at Devlin's proof of the exact big Ramsey degrees of the rationals,
as it has bearing on many current approaches to big Ramsey degrees.

The rationals can be represented by  the tree  $2^{<\om}$ of binary sequences with the
lexicographic order $\trl$ defined as follows:
Given $s,t\in 2^{<\om}$ with $s\ne t$, and   letting $u$ denote $s\wedge t$, define
 $s\trl t$  to  hold if and only if
($|u|<|s|$ and $s(|u|)=0$) or
$(|u|<|t|$ and $t(|u|)=1$).
Then  $(2^{<\om},\trl)$  is  a dense linear order.
The following is Definition 6.11 in \cite{TodorcevicBK10}, using the terminology of \cite{Sauer06}.
For $|s|<|t|$, the number $t(|s|)$ is called the {\em passing number} of $t$ at $s$.

\begin{definition}\label{defn.strongsimtype}
For $A,B\sse \om^{<\om}$,
we say that $A$ and $B$  are {\em similar} if there is a bijection $f:\mc(A)\ra \mc(B)$
such that for all $s,t\in \mc(A)$,
\begin{enumerate}
\item[(a)]  (preserves end-extension)
$s\sqsubseteq t \Leftrightarrow f(s)\sqsubseteq f(t)$,
\item[(b)] (preserves relative lengths)
$|s|<|t| \Leftrightarrow |f(s)|<|f(t)|$,
\item[(c)]
$s\in A \Leftrightarrow  f(s)\in B$,
\item[(d)] (preserves passing numbers)
$t(|s|)=f(t)(|f(s)|)$ whenever $|s|<|t|$.
\end{enumerate}
\end{definition}

Similarity is an equivalence relation;
a similarity equivalence class  is called  a {\em similarity type}.
We now outline the four  main steps to Devlin's characterization of big Ramsey degrees in the rationals.
Fix $n\ge 1$.

I. (Envelopes)
Given a  subset  $A\sse 2^{<\om}$ of size $n$,
let $k$ be the number of levels in $\mc(A)$.
An {\em envelope} of $A$ is a $k$-strong subtree  $E(A)$ of $2^{<\om}$ such that
$A\sse E(A)$.
Given any $k$-strong subtree $S$ of $2^{<\om}$, there is exactly one subset  $B\sse S$ which is similar to $A$.
This makes it possible to transfer a coloring of the similarity copies of $A$ in $2^{<\om}$
to the $k$-strong subtrees of $2^{<\om}$ in a well-defined manner.

II. (Finite Big Ramsey Degrees)
Apply  Milliken's theorem to obtain an infinite strong subtree $T\sse 2^{<\om}$ such that
every similarity copy of $A$ in $T$ has the same color.
As there are only finitely many similarity types of sets of size $n$,
finitely many applications of Milliken's theorem results in an infinite strong subtree $S\sse 2^{<\om}$  such that the coloring is monochromatic on each similarity type of  size $n$.
This achieves finite big Ramsey degrees.

III. (Diagonal Antichain  for Better Upper Bounds)
 To obtain the exact big Ramsey degrees,
 Devlin constructed a particular
antichain of nodes $D\sse 2^{<\om}$  such that
$(D,\trl)$ is a dense linear order and
no two nodes in
 the meet closure of $D$ have  the same length, a property called {\em diagonal}.
He also required  $(*)$: All passing numbers at the level of a terminal  node or a  meet node in $\mc(D)$ are $0$, except of course the rightmost extension of the meet node.
 Diagonal antichains turn out to be essential  to  characterizing big Ramsey degrees, whereas the additional requirement ($*$) is now seen to be
non-essential   when viewed through the lens of  coding trees of $1$-types.

IV.  (Exact Big Ramsey Degrees)
To characterize the big Ramsey degrees,
Devlin proved  that
the similarity type of each   subset of  $D$  of size $n$
persists in every  subset $D'\sse D$  such that  $(D',\trl)$ is a dense linear order.
The similarity types of antichains in $D$ thus form a   canonical partition for linear orders of size $n$.
By calculating the number of different similarity types of subsets of $D$ of size $n$, Devlin found the  big Ramsey degrees for the rationals.

Now we  present the characterization of the  big Ramsey degrees  for the rationals using coding trees of $1$-types.
Coding trees on $2^{<\om}$ were first developed in  \cite{DobrinenJML20}
to solve the  problem of whether or not the  triangle-free homogeneous graph has finite big Ramsey degrees.
The presentation given here is from \cite{CDP21}, where the notion of coding trees  was honed using model-theoretic ideas.
We hope that presenting this view here will set the stage for a
concrete understanding of
big Ramsey degree characterizations  discussed  in Section \ref{sec.RecentResults}.

Fix an enumeration   $\{q_0,q_1,\dots\}$ of $\mathbb{Q}$.
For $n< \om$, we let $\mathbb{Q}\re n$ denote the  substructure
$(\{q_i:i\in n\},<)$
of $(\bQ,<)$,
which
we refer to  as an {\em initial substructure}.
One can think of $\mathbb{Q}\re n$ as a finite approximation in a construction of the rationals.
The definition of a coding tree of $1$-types in \cite{CDP21} uses complete realizable quantifier-free $1$-types  over initial substructures.
Here, we shall retain the terminology of
\cite{CDP21}
 but (with apologies to model-theorists) will use sets of literals instead,
since this will  convey the important aspects of the constructions while being more accessible to a general readership.
For now, we  call a set of formulas $s\sse\{(q_i< x):i\in n\}\cup\{(x<  q_i):i\in n\}$   a {\em  $1$-type} over $\bQ\re n$ if
(a)
for each $i<n$ exactly one of the formulas
 $(q_i< x)$ or $(x<q_i)$ is in $s$,
and
(b)
  there is
some  (and hence infinitely many)
$j\ge n$ such that  $q_j$ {\em satisfies}  $s$,
meaning that
 replacing the variable $x$ by the rational number  $q_j$ in each  formula in $s$ results in  a true statement.
In other words,  $s$ is a $1$-type if $s$
prescribes  a  legitimate way to extend  $\bQ\re n$ to a linear order of size $n+1$.

\begin{definition}[Coding  Tree of $1$-Types for $\mathbb{Q}$, \cite{CDP21}]\label{defn.treecodeK}
For a fixed  enumeration $\{q_0,q_1,\dots\}$  of the rationals,
the {\em coding  tree of $1$-types}
$\bS(\mathbb{Q})$
 is the set of all
  $1$-types over initial substructures
along with a function $c:\om\ra \bS(\mathbb{Q})$ such that
$c(n)$ is the
$1$-type
of $q_n$
over $\mathbb{Q}\re n$.
The tree-ordering is simply inclusion.
\end{definition}

%%%%%%%%%%%%%%%
%%%%%%%%%%%%%%%

\begin{figure}[t]
\begin{tikzpicture}[grow'=up,scale=.6]
\tikzstyle{level 1}=[sibling distance=4in]
\tikzstyle{level 2}=[sibling distance=2in]
\tikzstyle{level 3}=[sibling distance=1in]
\tikzstyle{level 4}=[sibling distance=0.5in]
\tikzstyle{level 5}=[sibling distance=0.2in]
\tikzstyle{level 6}=[sibling distance=0.1in]
\tikzstyle{level 7}=[sibling distance=0.07in]
\node [label=$c_0$] {} coordinate (t9)
child{ coordinate (t0) edge from parent[color=black,thick]
child{ coordinate (t00) edge from parent[color=black,thick]
%edge from parent node[right] {help}
child{ coordinate (t000) edge from parent[color=black,thick]
child{coordinate (t0000) edge from parent[color=black,thick]
child{coordinate (t00000) edge from parent[color=black,thick]
child{coordinate (t000000) edge from parent[color=black,thick]}
}%endparen for t00000
}%endparen for t0000
}%endparen for t000
child{ coordinate (t001) edge from parent[color=black,thick]
child{ coordinate (t0010) edge from parent[color=black,thick]
child{coordinate (t00100) edge from parent[color=black,thick]
child{coordinate (t001000) edge from parent[color=black,thick]}
child{coordinate (t001001) edge from parent[color=black,thick]}
}%endparen for t00100
}%endparen for t0010
}%closeparen for t001
}%endparen for t00
%child{coordinate (t01) edge from parent[color=white]}
}%endparen for t0
child{ coordinate (t1) edge from parent[color=black,thick]
child{ coordinate (t10) edge from parent[color=black,thick]
child{ coordinate (t100)
child{coordinate (t1000)
child{coordinate (t10000)
child{coordinate (t100000)}
}%endparen for t10000
}%endparen for t1000
child{coordinate (t1001)
child{coordinate (t10010)
child{coordinate (t100100)}
}%endparen for t10010
}%endparen for t1001
}%endparen for t100
}%endparen for t10
child{ coordinate (t11) [label={\small 1/16}]
child{coordinate (t110)
child{coordinate (t1100)
child{coordinate (t11000)
child{coordinate (t110000)}
}%endparen t11000
child{coordinate (t11001)
child{coordinate (t110010)}
}%endparen for t11001
}%endparen for t1100
}%endparen for t110
}%endparen for t11
}%endparen for t1
;
%line spaces make a difference!!!

\node[circle, fill=blue,inner sep=0pt, minimum size=5pt] at (t9) {};
\node[circle, fill=blue,inner sep=0pt, minimum size=5pt,label=0:$c_2$,label=250:$\scriptstyle{x<q_1}$] at (t00) {};
\node[circle, fill=blue,inner sep=0pt, minimum size=5pt,label=0:$c_5$,label=280:${\scriptstyle x<q_4}$] at (t00100) {};
\node[circle, fill=blue,inner sep=0pt, minimum size=5pt, label=$c_1$,label=250:$\scriptstyle{q_0<x}$] at (t1) {};
\node[circle, fill=blue,inner sep=0pt, minimum size=5pt, label = 0:$c_3$,label=240:$\scriptstyle{q_2<x}$] at (t100) {};
\node[circle, fill=blue,inner sep=0pt, minimum size=5pt,label = 0:$c_4$,label=300:$\scriptstyle{q_3<x}$] at (t1100) {};
\node[label=280:$\scriptstyle{x< q_0}$] at (t0) {};
\node[label=270:$\scriptstyle{q_1<  x}$] at (t11) {};
\node[label=270:$\scriptstyle{x<  q_1}$] at (t10) {};
\node[label=300:$\scriptstyle{q_2 <  x}$] at (t110) {};
\node[label=240:${\scriptstyle x< q_3}$] at (t1000) {};
\node[label=280:${\scriptstyle q_3<  x}$] at (t1001) {};
\node[label=240:${\scriptstyle x< q_4}$] at (t10000) {};
\node[label=300:${\scriptstyle x< q_4}$] at (t10010) {};
\node[label=240:${\scriptstyle q_5<x}$] at (t100000) {};
\node[label=300:${\scriptstyle q_5< x}$] at (t100100) {};
\node[label=260:${\scriptstyle x<q_4}$] at (t11000) {};
\node[label=280:${\scriptstyle q_4<  x}$] at (t11001) {};
\node[label=260:${\scriptstyle q_5< x}$] at (t110000) {};
\node[label=280:${\scriptstyle q_5< x}$] at (t110010) {};
\node[label=260:${\scriptstyle x< q_2}$] at (t000) {};
\node[label=280:${\scriptstyle q_2< x}$] at (t001) {};
\node[label=260:${\scriptstyle x<q_3}$] at (t0000) {};
\node[label=280:${\scriptstyle x< q_3}$] at (t0010) {};
\node[label=260:${\scriptstyle x<q_4}$] at (t00000) {};
\node[label=260:${\scriptstyle x< q_5}$] at (t000000) {};
\node[label=260:${\scriptstyle x<q_5}$] at (t001000) {};
\node[label=280:${\scriptstyle q_5< x}$] at (t001001) {};
\node[label=280:${\scriptstyle q_5< x}$] at (t001001) {};

%the dots going up
%\node[circle,inner sep=0pt, minimum size=5pt,label=90:$\vdots$] at (t000000) {};
%\node[circle,inner sep=0pt, minimum size=5pt,label=90:$\vdots$] at (t001000) {};
%\node[circle,inner sep=0pt, minimum size=5pt,label=90:$\vdots$] at (t001001) {};
%\node[circle,inner sep=0pt, minimum size=5pt,label=90:$\vdots$] at (t100000) {};
%\node[circle,inner sep=0pt, minimum size=5pt,label=90:$\vdots$] at (t100100) {};
%\node[circle,inner sep=0pt, minimum size=5pt,label=90:$\vdots$] at (t100100) {};
%\node[circle,inner sep=0pt, minimum size=5pt,label=90:$\vdots$] at (t110000) {};
%\node[circle,inner sep=0pt, minimum size=5pt,label=90:$\vdots$] at (t110010) {};%#current

%empty circle nodes
\node[circle, fill=white,draw,inner sep=0pt, minimum size=4pt] at (t0) {};
\node[circle, fill=white,draw,inner sep=0pt, minimum size=4pt] at (t000) {};
\node[circle, fill=white,draw,inner sep=0pt, minimum size=4pt] at (t0000) {};
\node[circle, fill=white,draw,inner sep=0pt, minimum size=4pt] at (t00000) {};
\node[circle, fill=white,draw,inner sep=0pt, minimum size=4pt] at (t000000) {};
\node[circle, fill=white,draw,inner sep=0pt, minimum size=4pt] at (t001) {};
\node[circle, fill=white,draw,inner sep=0pt, minimum size=4pt] at (t0010) {};
\node[circle, fill=white,draw,inner sep=0pt, minimum size=4pt] at (t001000) {};
\node[circle, fill=white,draw,inner sep=0pt, minimum size=4pt] at (t001001) {};
\node[circle, fill=white,draw,inner sep=0pt, minimum size=4pt] at (t10) {};
\node[circle, fill=white,draw,inner sep=0pt, minimum size=4pt] at (t1000) {};
\node[circle, fill=white,draw,inner sep=0pt, minimum size=4pt] at (t10000) {};
\node[circle, fill=white,draw,inner sep=0pt, minimum size=4pt] at (t100000) {};
\node[circle, fill=white,draw,inner sep=0pt, minimum size=4pt] at (t1001) {};
\node[circle, fill=white,draw,inner sep=0pt, minimum size=4pt] at (t10010) {};
\node[circle, fill=white,draw,inner sep=0pt, minimum size=4pt] at (t100100) {};
\node[circle, fill=white,draw,inner sep=0pt, minimum size=4pt] at (t11) {};
\node[circle, fill=white,draw,inner sep=0pt, minimum size=4pt] at (t110) {};
\node[circle, fill=white,draw,inner sep=0pt, minimum size=4pt] at (t11000) {};
\node[circle, fill=white,draw,inner sep=0pt, minimum size=4pt] at (t110000) {};
\node[circle, fill=white,draw,inner sep=0pt, minimum size=4pt] at (t11001) {};
\node[circle, fill=white,draw,inner sep=0pt, minimum size=4pt] at (t110010) {};

%the nodes along the bottom
\node[circle, fill=blue,inner sep=0pt, minimum size=5pt, label=$q_2$,below=2.5cm of t00] (v2) {};
\node[circle, fill=blue,inner sep=0pt, minimum size=5pt,label=$q_5$,below =3.4cm of t001] (v5) {};
\node[circle, fill=blue,inner sep=0pt, minimum size=5pt,label=$q_0$,below =.7cm of t9] (v0) {};
\node[circle, fill=blue,inner sep=0pt, minimum size=5pt,label=$q_3$,below =3.4cm of t100] (v3) {};
\node[circle, fill=blue,inner sep=0pt, minimum size=5pt,label=$q_1$,below =1.6cm of t1] (v1) {};
\node[circle, fill=blue,inner sep=0pt, minimum size=5pt,label=$q_4$,below =2.5cm of t11] (v4) {};
\end{tikzpicture}
\caption{Coding tree $\bS(\bQ)$ of $1$-types for $(\bQ,< )$ and the  linear order represented by its coding nodes.}\label{fig.Qtree}
\end{figure}

%%%%%%%%%%%%%%%%%%%
%%%%%%%%%%%%%%%%%%%

Given $s\in \bS(\mathbb{Q})$
let $|s|=j+1$ where $j$ is maximal  such that
one of $(x< q_j)$ or $(q_j< x)$ is in $s$.
For each $i<|s|$,
  we let
$s(i)$ denote  the
 formula from among $(x< q_i)$ or $(q_i< x)$ which is in $s$.
%In this way, one can view $s$ as the sequence  $\lgl s(0),\dots, s(|s|-1)\rgl$.
The {\em coding nodes}  $c(n)$, in practice usually denoted by $c_n$, are special distinguished nodes representing the rational numbers;
$c_n$ represents the rational $q_n$, because $c_n$ is  the  $1$-type with  parameters
 from among $\{q_i:i\in n\}$ that $q_n$ satisfies.
Notice that this tree  $\bS(\mathbb{Q})$ has at most one splitting node per level.
The effect is that any antichain of coding nodes
in $\bS(\mathbb{Q})$ will automatically be diagonal.
(See Figure 1, reproduced from \cite{CDP21}.)

Fix an ordering  $<_{\mathrm{lex}}$ on the literals:
For $i<j$,  define
$(x< q_i)<_{\mathrm{lex}} (q_i< x)<_{\mathrm{lex}} (x< q_j)$.
Extend $<_{\mathrm{lex}}$ to
$\bS(\mathbb{Q})$ by declaring
for $s,t\in \bS(\mathbb{Q})$,
$s<_{\mathrm{lex}} t$ if and only if  $s$ and $t$ are incomparable
and for $i=|s\wedge t|$,
$s(i)<_{\mathrm{lex}} t(i)$.

\begin{definition}\label{def.simtypeQtypetree}
For $A,B$ sets of coding nodes in  $\bS(\mathbb{Q})$,
we say that $A$ and $B$ are {\em similar} if there is a bijection $f:\mc(A)\ra\mc(B)$ such that  for all $s,t\in\mc(A)$,
$f$ satisfies  (a)--(c) of Definition \ref{defn.strongsimtype}
and
(d${}'$)
$s <_{\mathrm{lex}}t \Leftrightarrow f(s) <_{\mathrm{lex}}f(t)$,
\end{definition}

When $B$ is similar to $A$, we call $B$ a
 {\em similarity copy} of $A$.
Condition (d) in  Definition \ref{defn.strongsimtype}
implies that the lexicographic order on $2^{<\om}$ is preserved,
and moreover,
 that passing numbers at
meet nodes and at terminal nodes   are preserved.
In (d${}'$)
we only need to preserve lexicographic order.

Extending Harrington's method,
 forcing is utilized to
obtain a pigeonhole principle for coding trees of $1$-types in the vein of the Halpern--\Lauchli\ Theorem \ref{thm.HL}, but for colorings of finite sets of coding nodes, rather than antichains.
Via an inductive argument using this pigeonhole principle, we  obtain the following Ramsey theorem on coding trees.

\begin{theorem}[\cite{CDP21}]\label{thm.SDAPQ}
Let  $\bS(\bQ)$ be
 a coding tree of $1$-types for the rationals.
Given a finite set $A$ of coding nodes in $\bS(\bQ)$ and a finite coloring of all similarity copies of $A$ in  $\bS(\bQ)$,
there is a  coding subtree $S$ of $\bS(\bQ)$ similar to $\bS(\bQ)$
 such that
all
similarity copies of $A$ in $S$  have the same color.
\end{theorem}

%is there a way to say this without saying similarity copy all the time?

Fix $n\ge 1$.
By applying Theorem \ref{thm.SDAPQ} once for each similarity type of coding nodes of size $n$,
we prove  finite  big Ramsey degrees, accomplishing step II while bypassing step I in Devlin's proof.
Upon taking {\em any}
 antichain $D$ of coding nodes in $\bS(\mathbb{Q})$
representing a dense linear order,
we obtain better upper bounds which are then proved to be exact, accomplishing steps III and IV.
\vskip.1in

\noindent \textbf{Big Ramsey degrees of the rationals}.
In \cite{CDP21}, we show that
given $n\ge 1$, the big Ramsey degree $T(n)$ for linear orders of size $n$ in the rationals is
the number of
similarity types of antichains of coding nodes in $\bS(\bQ)$.

What then {\em is} the big Ramsey degree  $T(n)$ in the rationals?
It  is the
number of different ways to order the indexes
of  an increasing sequence of  rationals
 $\{q_{i_0}< q_{i_1}<\dots< q_{i_{n-1}}\}$
 with
 incomparable $1$-types
along with the  number of ways to order the
 first differences  of their $1$-types
over initial substructures of $\mathbb{Q}$.
The first difference    between the $1$-types of  the rationals  $q_i$ and $q_j$
occurs at  the least $k$ such that $q_i < q_k $ and $q_k<  q_j$, or vice versa.
This means that
$q_i$ and $q_j$
are in
the same interval of  $\bQ\re k$
but in  different intervals of $\bQ\re (k+1)$.
Concretely, $T(n)$ is the number of
$<$-isomorphism classes of $(2n-1)$-tuples   of integers
$(i_0,\dots,i_{n-1},k_0,\dots,k_{n-2})$
with the following properties:
 $\{q_{i_0}< q_{i_1}< \dots< q_{i_{n-1}}\}$ is a set of rationals in increasing order, and  for each $j<n-1$,
$q_{i_j}< q_{k_j}< q_{i_{j+1}}$  where  $k_j<\min(i_j,i_{j+1})$ and  is the least  integer  satisfying this relation.

%%%%%%%%%%%%%%%%
%%%%%%%%%%%%%%%%
%%%%%%%%%%%%%%%%
%%%%%%%%%%%%%%%%

\section{Historical highlights, recent results, and methods}\label{sec.HH}

We now  highlight some historical achievements,
and present recent results and the main ideas of their methods.
For an   overview of results up to the year 2000, see  the appendix  by Sauer in
 \Fraisse's book   \cite{FraisseBK};
for an overview up to the year 2013, see Nguyen Van Th\'{e}'s
habilitation thesis \cite{NVTHabil}.
Those interested in open problems intended for undergraduate research may enjoy
 \cite{DG20}.

The Rado graph is the second  example of
a homogeneous structure
with non-trivial big Ramsey degrees
 which has been fully understood in terms of its partition theory.
The Rado graph $\bR$
is up to isomorphism  the homogeneous graph on countably many vertices which is universal for all countable graphs.
It was known to   \Erdos\ and other Hungarian mathematicians in the 1960's, though possibly earlier,
that the
the Rado graph is indivisible.
In
their
1975 paper  \cite{EHP73},
\Erdos, Hajnal, and P\'{o}sa
constructed a  coloring of the edges in $\bR$ into two colors such that
both colors persist
in each  subcopy of $\bR$.
Pouzet and Sauer later showed
in \cite{Pouzet/Sauer96}
 that
 the big Ramsey degree for  edge colorings
 in the Rado graph  is exactly two.
The complete characterization of the big Ramsey degrees of the Rado graph was achieved  in a pair of papers by
 Sauer \cite{Sauer06} and by Laflamme, Sauer, and Vuksanovic \cite{LSV06}, both appearing in 2006,
and the degrees were calculated by Larson in \cite{Larson08}.
The two papers
\cite{Sauer06} and \cite{LSV06}  in fact characterized  exact big Ramsey degrees
for all  unrestricted homogeneous structures
with finitely many binary relations, including   the homogeneous digraph, homogeneous tournament, and random graph with finitely many edges of different colors.
Milliken's Theorem was used to prove existence of  upper bounds,   alluding to a deep connection between big Ramsey degrees and Ramsey theorems for trees.  These  results are discussed  in Subsection \ref{subsec.SDAP}.

In \cite{LNVTS10},  for each $n\ge 2$,
Laflamme, Nguyen Van Th\'{e},  and Sauer
calculated the big Ramsey degrees of  $\mathbb{Q}_n$,
the rationals  with an equivalence relation with $n$ many equivalence classes each of which is dense in
$\mathbb{Q}$.
This hinged on  proving a ``colored version" of  Milliken's theorem,
where the levels of the trees are colored,  to achieve upper bounds.
Applying  their result for $\mathbb{Q}_2$, they calculated  the big Ramsey degrees of   the
dense local order, denoted  $\bfS(2)$.
In  his PhD thesis \cite{HoweThesis},
Howe   proved
finite big Ramsey degrees for the generic bipartite graph and the \Fraisse\ limit of the class of finite linear orders with a convex equivalence relation.

A robust and streamlined approach
applicable to  a large class of homogeneous structures, and
recovering the previously mentioned examples
(except for
$\mathbf{S}(2)$),
was developed
by Coulson, Patel, and the author
in
\cite{CDP21}, building on ideas in \cite{DobrinenJML20} and \cite{DobrinenH_k19}.
In \cite{CDP21}, it was shown that
homogeneous structures with relations of arity at most two satisfying  a strengthening of SAP, called SDAP$^+$, have big Ramsey structures  which are characterized in a simple manner, and therefore
their big Ramsey degrees  are easy to compute.
The proof  proceeds
via a Ramsey theorem for colorings of finite antichains of coding nodes
on
 {\em diagonal coding trees of $1$-types}.
This approach bypasses any need for envelopes, the  theorem producing of its own accord exact upper bounds.
Moreover, the
Halpern--\Lauchli-style theorem, which is proved  via
 forcing arguments to achieve a ZFC result and used as the
 pigeonhole principle
 in the Ramsey  theorem,
 immediately yields
indivisibility for   all
 homogeneous structures  satisfying
SDAP$^+$,
with relations of any arity.
These results and their methods are discussed in
Subsection \ref{subsec.SDAP}.

The $k$-clique-free  homogeneous graphs, denoted $\mathbf{G}_k$, $k\ge 3$,
 were constructed by Henson in  his 1971 paper \cite{Henson71}, where he proved these graphs to be weakly indivisible.
In their 1986 paper \cite{Komjath/Rodl86},
Komj\'{a}th and \Rodl\
proved that $\mathbf{G}_3$ is indivisible,  answering a question of Hajnal.
A few years later,
 El-Zahar and Sauer gave a systematic approach  in  \cite{El-Zahar/Sauer89},
proving that for each $k\ge 3$, the  $k$-clique-free homogeneous graph $\mathbf{G}_k$ is indivisible.
In 1998,
Sauer proved in
\cite{Sauer98} that  the big Ramsey degree for edges in   $\mathbf{G}_3$
is two.
Further  progress on big Ramsey degrees of $\mathbf{G}_3$, however, needed a new approach.
This was achieved by the author in
 \cite{DobrinenJML20},
where  the   method of coding trees was first developed.
In \cite{DobrinenH_k19}, the author
extended this work, proving
that $\mathbf{G}_k$ has finite big Ramsey degrees, for each $k\ge 3$.
In  \cite{DobrinenJML20} and \cite{DobrinenH_k19},  the author proved a Ramsey theorem
for colorings of finite antichains of  coding nodes
in  diagonal coding trees.
These diagonal coding trees were designed to
achieve very good upper bounds
and directly recover the indivisibility results in
 \cite{Komjath/Rodl86} and \cite{El-Zahar/Sauer89},
 discovering much of the essential structure involved in characterizing their exact big Ramsey degrees.
(Milliken-style theorems on non-diagonal coding trees  which fully branch at each level do not directly prove indivisibility results, and produce looser upper bounds.)
In particular,
after a minor  modification,  the trees
in  \cite{DobrinenJML20}
 produced exact big Ramsey degrees for $\mathbf{G}_3$, as shown in
 \cite{DobrinenH_3Exact20}.
Around the same time,
exact big Ramsey degrees for  $\mathbf{G}_3$ were    independently proved  by
Balko,  Chodounsk{\'{y}},
 Hubi{\v{c}}ka,  Kone{\v{c}}n{\'{y}},  Vena,  and Zucker, instigating the  collaboration of this group with the author.

Given a finite relational language $\mathcal{L}$,
an $\mathcal{L}$-structure $\bfA$ is called {\em irreducible} if
each pair of its vertices are in some relation of $\bfA$.
Given a set $\mathcal{F}$ of finite irreducible $\mathcal{L}$-structures,
Forb$(\mathcal{F})$ denotes
the class of all finite $\mathcal{L}$-structures into which  no member of $\mathcal{F}$ embeds.
\Fraisse\ classes
of the form Forb$(\mathcal{F})$
are exactly those with free amalgamation.
Zucker in
\cite{Zucker20}  proved that
for any \Fraisse\ class of the form  Forb$(\mathcal{F})$,
where $\mathcal{F}$ is  a finite set of irreducible substructures and all
 relations have arity at most two,
its \Fraisse\ limit has
 finite big Ramsey degrees.
His proof used coding trees which branch at each level and  a forcing argument to obtain a
Halpern--\Lauchli-style theorem which formed the pigeonhole principle for a  Milliken-esque theorem for these coding trees.
An important advance  in this paper is Zucker's abstract, top-down  approach,
providing  simplified and relatively short proof of finite big Ramsey degrees for
this large class of homogeneous structures.
On the other hand,
 his Milliken-style theorem   does not directly
 recover indivisibility (more work is needed afterwards to show this), and the upper bounds in
\cite{Zucker20}  did not recover those in
\cite{DobrinenJML20} or
\cite{DobrinenH_k19} for the  homogeneous $k$-clique-free graphs.
However, by
 further work done
in
\cite{Balko7_binaryFAP},
by
Balko,  Chodounsk{\'{y}},
 Hubi{\v{c}}ka,  Kone{\v{c}}n{\'{y}},  Vena,  Zucker, and the author,
indivisibility results are proved  and
 exact big Ramsey degrees  are characterized.
Thus, the picture for FAP classes with finitely many  relations of arity at most two is now clear.
These results will be discussed in Subsection \ref{subsec.FAP}.

Next, we look at   homogeneous structures with relations of arity at most two  which do not satisfy SDAP$^+$ and  whose ages have strong (but not free) amalgamation.
 Nguyen Van Th\'{e}
made a  significant contribution
  in  his 2008 paper
\cite{NVT08}, in which he proved  that the
 ultrametric Urysohn space  $\mathbf{Q}_S$ has finite big Ramsey degrees if and only if $S$ is a finite distance set.
In the case that $S$ is  finite, he calculated the big Ramsey degrees.
Moreover, he showed that for an infinite countable distance set $S$,  $\mathbf{Q}_S$ is indivisible if and only if $S$ with the reverse order as a subset of the reals is well-ordered.
His proof used infinitely wide trees of finite height and his  pigeonhole principle was actually Ramsey's theorem.
All countable Urysohn metric spaces with  finite distance set
were proved to be
 indivisible by
Sauer in \cite{Sauer12},
completing the work that was initiated in  \cite{NVTS09}
in relation to the celebrated distortion problem from Banach space theory and its solution by Odell and Schlumprecht in \cite{OS94}.

%I MAY WANT TO ADD MORE NOW THAT I KNOW MORE HERE.

  \Masulovic\  instigated the use of category theory to prove transport principles
showing that finite big Ramsey degrees can be inferred from one category to another.
After proving a general
 transport principle
 in  \cite{Masulovic18}, he  applied it to prove finite big Ramsey degrees for many universal structures and also
for homogenous metric spaces with
 finite distance sets  with a certain property which he calls {\em compact with one nontrivial block}.
\Masulovic\
proved in \cite{Masulovic_RBS20}  that in categories satisfying certain
 mild conditions,   small Ramsey degrees
are minima of big Ramsey degrees.
In the  paper  \cite{MS} with \v{S}obot (not using category theory),
finite big Ramsey degrees for finite chains in countable ordinals
were shown to  exist if and only if the ordinal is smaller than $\om^{\om}$.
Dasilva Barbosa in \cite{Barbosa20}
proved that categorical precompact expansions grant upper bounds for big and small Ramsey degrees.
As an application, he calculated the big Ramsey degrees of the circular directed graphs $\mathbf{S}(n)$ for all $n\ge 2$,
extending the work  in \cite{LNVTS10}  for  $\mathbf{S}(2)$.

 \Hubicka\ recently developed a new method  to handle forbidden substructures utilizing topological Ramsey spaces of parameter words due to Carlson and Simpson \cite{Carlson/Simpson84}.
In  \cite{Hubicka_CS20},  he applied his method
to prove that the
homogeneous partial order and
Urysohn $S$-metric spaces
(where $S$ is a set of non-negative reals with $0\in S$ satisfying the 4-values condition)
have finite big Ramsey degrees.
He also showed that this method
is quite broad and
 can be applied to yield
 a short proof of  finite big Ramsey degrees  in $\mathbf{G}_3$.
%We note that bounds in \cite{DobrinenJML20} were not immediately recovered by this method.
Beginning with the upper bounds in
 \cite{Hubicka_CS20},
 the exact big Ramsey degrees of
 the generic partial order
have been  characterized in
\cite{Balko7_PO} by
Balko,  Chodounsk{\'{y}},
 Hubi{\v{c}}ka,  Kone{\v{c}}n{\'{y}},  Vena,   Zucker, and the author.
Also utilizing techniques from
 \cite{Hubicka_CS20},
Balko, Chodounsk{\'{y}},
Hubi{\v{c}}ka,  Kone{\v{c}}n{\'{y}}, Ne\v{s}et\v{r}il, and  Vena in \cite{Balko6_ForbCycle}
have  found a condition which guarantees finite big Ramsey degrees for  binary relational homogeneous structures  with strong amalgamation.
Examples of structures satisfying this condition include the $S$-Urysohn space  for finite distance sets $S$, $\Lambda$-ultrametric spaces for a finite distributive lattice,
and metric spaces associated to metrically homogeneous graphs of a finite diameter from Cherlin's list with no Henson constraints.

For homogeneous structures with free amalgamation  with finitely many  relations of any  arity, the picture for indivisibility  is now  clear due  to the recent breakthrough  of Sauer.
In  \cite{Sauer20}, Sauer showed that a homogeneous free amalgamation structure with relations of finite arity is indivisible if and only if its age poset is linearly ordered, a property he called {\em rank linear},
 culminating  a line of work  in  \cite{El-Zahar/Sauer91}, \cite{El-Zahar/Sauer93},  \cite{El-Zahar/Sauer94},  \cite{Sauer03}, and  \cite{El-Zahar/Sauer05}.
On the other hand,
big Ramsey degrees of structures with relations of arity greater than two
has only recently
seen progress,
 beginning with
 \cite{Balko5_3unifGraphs19}, where
Balko, Chodounsk{\'{y}},
Hubi{\v{c}}ka,  Kone{\v{c}}n{\'{y}},  and Vena
found upper bounds for the big Ramsey degrees of
 the generic $3$-hypergraph.
Work in this area is ongoing and promising.

%%%%%%%%%%%%%%%%
%%%%%%%%%%%%%%%%
%%%%%%%%%%%%%%%%
%%%%%%%%%%%%%%%%

\section{Exact big Ramsey degrees}\label{sec.RecentResults}

This section presents
characterizations  of
exact big Ramsey degrees
 known at the time of writing.
These hold
for
 homogeneous structures with finitely many relations of arity at most two.
Two general classes have been completely understood:
Structures satisfying a certain strengthening of strong amalgamation  called SDAP$^+$ (Subsection \ref{subsec.SDAP}) and
structures whose ages have free amalgamation (Subsection \ref{subsec.FAP}).
Lying outside of these two classes,
the generic partial order has been completely understood in terms of exact big Ramsey degrees and will  be   briefly discussed  at the end of
Subsection \ref{subsec.FAP}.
These characterizations  all involve the notion of a diagonal  antichain,
in  various  trees or spaces of parameter words,
representing a copy of an enumerated homogeneous structure.
Here, we present these
notions
in terms of structures, as
they
are independent of the representation.

%\begin{definition}\label{defn.antichaindiag}
Let $\bK$ be an  enumerated homogeneous structure with universe $\{v_n:n<\om\}$.
Let $\bfA\le\bK$ be a finite substructure of $\bK$,
and suppose that the universe of $\bfA$ is
$\{v_i:i\in I\}$ for some finite set $I\sse \om$.
We say that $\bfA$ is an {\em antichain} if for each pair $i<j$ in $I$ there is a $k(i,j)<i$ such that  the set
 $\{k(i,j):i,j\in I\mathrm{\ and\ }i<j\}$ is disjoint from $I$,  and
\begin{align}
\bK\re (\{v_{\ell}:\ell<k(i,j)\}\cup\{v_{i}\})
& \cong
\bK\re (\{v_{\ell}:\ell<k(i,j)\}\cup\{v_{j}\})\cr
\bK\re (\{v_{\ell}:\ell\le k(i,j)\}\cup\{v_{i}\})
& \not\cong
\bK\re (\{v_{\ell}:\ell\le k(i,j)\}\cup\{v_{j}\}).
\end{align}
An antichain $\bfA$ is called {\em diagonal}
if $\{k(i,j):i<j\le m\}$ has cardinality $m$.
We call $k(i,j)$ the {\em meet level} of the pair $v_{i},v_{j}$.
%\end{definition}

The notion of diagonal antichain is central to all  characterizations of big Ramsey degrees obtained so far.
It seems likely that antichains will be  essential to all characterizations of big Ramsey degrees.
However,
preliminary work  shows that some homogeneous binary relational structures, such as  two or more independent linear orders,
will    have  characterizations  in their trees of $1$-types
 involving antichains which are not diagonal,
but could still be characterized via products of finitely many diagonal antichains.

The indexing  of the relation symbols
$\{R_{\ell}:\ell<L\}$ in the
 language
$\mathcal{L}$ of $\bK$
 induces a lexicographic ordering on trees representing relational structures.
Here, we present this idea directly on the structures.
For $m\ne n$,
we declare $v_m<_{\mathrm{lex}}v_n$ if and only if
$\{v_m,v_n\}$ is an antichain and,
letting
$k$ be the meet level of the pair
$v_m,v_n$,
and letting
$\ell$ denote the least index in $L$ such that $v_m$ and $v_n$ disagree on their $R_{\ell}$-relationship with $v_{k}$,
either $R_{\ell}(v_k,v_n)$ holds while
$R_{\ell}(v_k,v_m)$ does not,
or else
$R_{\ell}(v_n,v_k)$ holds while
$R_{\ell}(v_m,v_k)$ does not.

Two diagonal antichains $\bfA$ and $\bfB$
in an enumerated homogeneous structure $\bK$
 are {\em similar}
if they have the same number of vertices,
and the increasing bijection from the universe $A=\{v_{m_i}:i\le p\}$  of $\bfA$ to the universe  $B=\{v_{n_i}:i\le p\}$ of $\bfB$ induces an
 isomorphism  from $\bfA$ to $\bfB$ which preserves  $<_{\mathrm{lex}}$
and  induces a map on  the meet levels which,
for each $i<j\le p$, sends $k(m_i,m_j)$ to $k(n_i,n_j)$.
This  implies  that  the
map  sending the coding node $c_{m_i}$ to $c_{n_i}$ ($i\le p$)
in the coding tree of $1$-types  $\bS(\bK)$
(see Definition \ref{defn.treecodeK})
induces a map
on the meet-closures
 of $\{c_{m_i}: i\le p\}$ and $\{c_{n_i}:i\le p\}$
satisfying  Definition
\ref{def.simtypeQtypetree}.

Similarity is an equivalence relation, and an equivalence class is called a {\em similarity type}.
We say that $\bK$ has {\em simply characterized big Ramsey degrees} if
for $\bfA\in\Age(\bK)$,
the big Ramsey degree of $\bfA$ is exactly
the number of similarity types of diagonal antichains representing  $\bfA$.
In the next Subsection, we will see many homogeneous structures with simply characterized big Ramsey degrees.

%%%%%%%%%%%%%%%%
%%%%%%%%%%%%%%%%

\subsection{Exact big Ramsey degrees with a simple characterization}\label{subsec.SDAP}

The decades-long investigation of the big Ramsey degrees of the Rado graph  culminated in
 the two papers  \cite{Sauer06} and \cite{LSV06}.
These two papers  moreover
characterized the big Ramsey degrees for all   {\em unrestricted} binary relational homogeneous structures.
Unrestricted
binary relational   structures
 are determined by a finite   language
$\mathcal{L}=\{R_0,\dots, R_{l-1}\}$ of binary relation symbols
and a non-empty constraint set
$\mathcal{C}$ of $\mathcal{L}$-structures with universe $\{0,1\}$
with the  following property:
 If $\bfA$ and $\bfB$ are two isomorphic $\mathcal{L}$-structures with universe $\{0,1\}$,
then either both are in $\mathcal{C}$ or neither is in
$\mathcal{C}$.
We let   $\mathbf{H}_{\mathcal{C}}$
denote  the homogeneous structure
such that  each of its   substructures
 with universe of size two is isomorphic to one of the structures in $\mathcal{C}$.
Examples of unrestricted binary relational  homogeneous structures include the Rado graph, the generic directed graph, the generic tournament,  and
random graphs with more than one edge relation.

Given  a universal constraint set  $\mathcal{C}$,
letting
 $k=|\mathcal{C}|$,
Sauer showed in \cite{Sauer06}
how to form a
structure, call it $\mathbf{U}_{\mathcal{C}}$, with nodes  in
 the  tree  $k^{<\om}$ as
vertices, such that
 $\mathbf{H}_{\mathcal{C}}$ embeds into
$\mathbf{U}_{\mathcal{C}}$.
Fix a bijection  $\lambda:\mathcal{C} \ra k$.
Given two nodes $s,t\in k^{<\om}$ with $|s|<|t|$,
declare that
$t(|s|)=j$
if and only if
the induced substructure of  $\mathbf{U}_{\mathcal{C}}$
on universe $\{s,t\}$
is  isomorphic to the structure $\lambda(j)$ in $\mathcal{C}$, where the isomorphism sends $s$ to $0$ and $t$ to $1$.
For two nodes $s,t\in k^{<\om}$ of the same length,  declare that for $s$ lexicographically less than $t$,
the induced substructure of  $\mathbf{U}_{\mathcal{C}}$
on universe $\{s,t\}$
is  isomorphic to the structure $\lambda(0)$ in $\mathcal{C}$, where the isomorphism sends $s$ to $0$ and $t$ to $1$.
As a special case,
 a universal graph is constructed as follows:
Let each node  in  $2^{<\om}$ be  a vertex.
Define an edge relation $E$  between vertices by
declaring that, for $s\ne t$ in
 $2^{<\om}$,
$s\, E\, t$ if and only if $|s|\ne |t|$ and
 ($|s|<|t| \Rightarrow  t(|s|)=1$).
Then $(2^{<\om},E)$ is universal for all countable graphs.
In particular, the Rado graph embeds into  the graph $(2^{<\om},E)$
 and vice versa.

In trees of the form $k^{<\om}$,
the notion of similarity is exactly that of Definition \ref{defn.strongsimtype},
 and
steps I--IV discussed in Section \ref{sec.CaseStudies}
outline the proof of exact big Ramsey degrees contained in the pair of papers \cite{Sauer06} and \cite{LSV06}.
Milliken's Theorem was  used to prove existence of  upper bounds via strong tree envelopes.
For step III,
 Sauer  constructed in \cite{Sauer06}  a diagonal antichain $D\sse k^{<\om}$
such that the substructure of $\mathbf{U}_{\mathcal{C}}$ restricted to universe $D$ is isomorphic to $\mathbf{H}_{\mathcal{C}}$, achieving
upper bounds shown to be exact in
 \cite{LSV06}, finishing step IV.
 The big Ramsey degree of a finite substructure $\bfA$  of $\mathbf{H}_{\mathcal{C}}$ is  exactly  the number of distinct similarity types of subsets of $D$ whose induced substructure in
$\mathbf{U}_{\mathcal{C}}$
is isomorphic to $\bfA$.

The work in \cite{Sauer06}   and  \cite{LSV06}
greatly influenced the author's
development of
 {\em coding trees} and their Ramsey theorems  in \cite{DobrinenJML20} and
\cite{DobrinenH_k19} (discussed in Subsection \ref{subsec.FAP}).
Those papers along with
a suggestion of Sauer  to the author during the Banff 2018 Workshop on {\em Unifying Themes in Ramsey Theory}, to try moving the forcing arguments in those papers from coding trees to structures,
informed the approach taken in  the paper
\cite{CDP21}, which is now discussed.

Let $\bK$ be an enumerated \Fraisse\ structure
with vertices $\{v_n:n<\om\}$.
For $n<\om$, we let $\bK_n$ denote $\bK\re \{v_i:i<n\}$, the  induced substructure of $\bK$ on its first $n$ vertices, and call $\bK_n$ an {\em initial substructure} of $\bK$.
We write {\em $1$-type} to mean complete realizable quantifier-free $1$-type over
$\bK_n$ for some  $n$.

\begin{definition}[Coding  Tree of $1$-Types, \cite{CDP21}]\label{defn.CT1Type}
The {\em coding  tree of $1$-types}
$\bS(\bK)$
for an enumerated \Fraisse\ structure $\bK$
 is the set of all
  $1$-types over initial substructures of $\bK$
along with a function $c:\om\ra \bS(\bK)$ such that
$c(n)$ is the
$1$-type
of $v_n$
over $\bK_n$.
The tree-ordering is simply inclusion.
\end{definition}

A substructure $\bfA$ of $\bK$ with universe $A=\{v_{n_0},\dots, v_{n_m}\}$ is represented by the set of coding nodes $\{c(n_0),\dots, c(n_m)\}$ as follows:
For each $i\le m$,
since $c(n_i)$ is the quantifier-free $1$-type
of $v_{n_i}$ over $\bK_{n_i}$,
substituting $v_{n_i}$ for the variable $x$ into each formula in $c(n_i)$ which has only  parameters  from $\{v_{n_j}:j<i\}$
uniquely determines the relations in $\bfA$ on  the vertices
$\{v_{n_j}:j\le i\}$.
In \cite{CDP21}, we formulated the following
 strengthening  of SAP in order to
extract  a  general property ensuring that
big Ramsey degrees have
 simple characterizations.

\begin{definition}[SDAP]\label{defn.EEAP_new}
A \Fraisse\ class $\mathcal{K}$ has the
{\em Substructure Disjoint Amalgamation Property (SDAP)} if $\mathcal{K}$
has strong amalgamation,
and the following holds:
Given   $\bfA, \bfC\in\mathcal{K}$,
 suppose that
 $\bfA$ is a substructure of $\bfC$, where   $\bfC$ extends  $\bfA$ by two vertices, say $v$ and $w$.
Then there exist  $\bfA',\bfC'\in\mathcal{K}$, where
$\bfA$ is a substructure of
$\bfA'$
and
$\bfC'$ is a disjoint amalgamation of $\bfA'$ and $\bfC$ over $\bfA$, such that
letting   $v',w'$ denote the two vertices in
 $C'\setminus A'$ and
assuming (1) and (2), the conclusion holds:
 \begin{enumerate}
 \item[(1)]
Suppose
$\bfB\in\mathcal{K}$  is any structure
 containing $\bfA'$ as a substructure,
and let
 $\sigma$ and $\tau$  be
  $1$-types over $\bfB$  satisfying    $\sigma\re\bfA'=\type(v'/\bfA')$ and $\tau\re\bfA'=\type(w'/\bfA')$,
\item[(2)]
Suppose
$\bfD\in \mathcal{K}$  extends  $\bfB$
by one vertex, say $v''$, such that $\type(v''/\bfB)=\sigma$.
\end{enumerate}
Then
  there is
an  $\bfE\in\mathcal{K}$ extending   $\bfD$ by one vertex, say $w''$, such that
  $\type(w''/\bfB)=\tau$ and  $\bfE\re (
  \mathrm{A}\cup\{v'',w''\})\cong \bfC$.
\end{definition}

This amalgamation property can of course be presented in terms of embeddings, but the form here is indicative of how it is utilized.
A free amalgamation version called
{\em SFAP} is obtained from SDAP by restricting to FAP classes and requiring
$\bfA'=\bfA$ and $\bfC'=\bfC$.
Both of these amalgamation properties are
preserved  under free superposition.
A {\em diagonal  subtree}  of $\bS(\bK)$  is a subtree such that at any level, at most one node branches, the branching degree is two,
and branching and coding nodes never occur on the same level.
Diagonal coding trees are subtrees of $\bS(\bK)$ which are diagonal
and represent a subcopy of $\bK$.
The property SDAP$^+$ holds for a homogeneous structure $\bK$ if
(a) its age satisfies SDAP,
(b) there is
a diagonal coding
 subtree of $\bS(\bK)$,
and
(c)  a technicality called the Extension Property which in most cases is trivially satisfied.
Classes of the form Forb$(\mathcal{F})$ where
$\mathcal{F}$ is a finite set of $3$-irreducible structures,
meaning each triple of vertices is in some relation,
satisfy SFAP;
their ordered versions satisfy SDAP$^+$.

A version of the Halpern--\Lauchli\ Theorem for  diagonal coding trees was proved in \cite{CDP21} using the method of forcing to obtain a ZFC result, with the following theorem as an immediate consequence.

\begin{theorem}[\cite{CDP21}]\label{thm.indivisibility}
Let $\bK$ be a homogeneous structure satisfying SDAP$^+$, with finitely many relations of any arity.
Then $\bK$ is indivisible.
\end{theorem}

For relations of arity at most two, an induction proof then yields a
Ramsey theorem for
 finite  colorings of finite antichains of coding nodes in  diagonal coding trees.
This
 accomplishes steps I--III simultaneously and directly,
without any need for envelopes, providing upper bounds which are then proved to be exact,
finishing step IV.

\begin{theorem}[\cite{CDP21}]\label{thm.CDPmain}
Let
$\bK$ be a homogeneous structure satisfying
SADP$^+$, with finitely many relations of arity at most two.
Then $\bK$
admits  a big Ramsey structure and moreover,
has simply characterized  big Ramsey degrees.
\end{theorem}

Theorem \ref{thm.CDPmain} provides new classes of examples
 of big Ramsey structures while
 recovering results in
 \cite{DevlinThesis},
\cite{HoweThesis},
 \cite{LNVTS10}, and
 \cite{LSV06}
and special cases of the results in \cite{Zucker20}.
Theorem \ref{thm.indivisibility} provides new  classes of examples of indivisible \Fraisse\ structures, in particular for ordered structures such as  the ordered Rado graph,
while recovering results in
\cite{El-Zahar/Sauer89},
\cite{Komjath/Rodl86},  and
\cite{El-Zahar/Sauer94} and some of the results in
 \cite{Sauer20}.
 Sauer's work in \cite{Sauer20}
provides the  full picture on indivisibility for
 free amalgamation classes.
 Theorem \ref{thm.indivisibility} recovers certain cases of Sauer's results for FAP classes,
while  providing new SAP examples  with  indivisibility.

%%%%%%%%%%%%%%%%%%%%%%%
%%%%%%%%%%%%%%%%%%%%%%%

\subsection{Big Ramsey degrees for free amalgamation classes}\label{subsec.FAP}

An obstacle to progress
in  partition theory of  homogeneous structures
 had been the fact that
Milliken's theorem is not able to handle forbidden substructures, for instance, triangle-free graphs.
Most results up to 2010 had either utilized
 Milliken's theorem or a variation
(as in \cite{Sauer06}, \cite{LNVTS10})
or else used difficult direct methods
 (as in \cite{Sauer98}) which did not lend naturally to  generalizations.
The idea of coding trees  came to the author
during the
her stay at the
 Isaac Newton Institute in 2015 for the
programme, {\em Mathematical, Foundational and Computational Aspects of the Higher Infinite},
culminating   in the work \cite{DobrinenJML20}.
The ideas behind coding trees
included the following:
Knowing that at the end of the process one will  want a diagonal antichain representing a copy of $\mathbf{G}_3$,
starting with a tree where  vertices in $\mathbf{G}_3$ are represented by special  nodes on  different levels should not hurt the results.
Further, by using special nodes to code the vertices of $\mathbf{G}_3$ into the trees,
one might have a chance to
prove
 Milliken-style theorems on  a
collection of trees, each of which
 codes a  subcopy of $\mathbf{G}_3$.

The author had made a previous   attempt at this problem starting early in 2012.
Upon  stating her  interest  this problem,
Todorcevic (2012, at the Fields Institute Thematic Program on Forcing and Its Applications)  and Sauer (2013, at the \Erdos\ Centenary Meeting) each told the author that a new kind of Milliken theorem would need to be developed in order to handle triangle-free graphs, which intrigued her even more.
Though unknown to her at the time, a key piece to this puzzle would be Harrington's forcing proof of the Halpern--\Lauchli\ theorem, which Laver was kind enough to outline to her in 2011.
(At that time, the author was unaware of the proof in \cite{Farah/TodorcevicBK}.)
While at the INI in 2015,
Barto\v{s}ov\'{a} reminded the author of her interest in big Ramsey degrees of $\mathbf{G}_3$.
Having had time by then  to fill out and digest Laver's outline,
it occurred to the author  to try
approaching the problem first with the strongest tool available, namely forcing.

Forcing is a set-theoretic method which is normally used to extend a given universe satisfying a given set of axioms (often ZFC)
to a larger universe in which  those same set of axioms hold while  some other statement or property is different than in  the original universe.
The beautiful thing about Harrington's proof is that, while it does involve the method of forcing,
the forcing is  only used as a search engine for an object which already  exists  in the universe  in  which one lives.
In the context of the \Fraisse\ limit  $\bK$ of
a class  $\mathrm{Forb}(\mathcal{F})$, where $\mathcal{F}$ is a finite set of finite irreducible structures,
by carefully designing forcings on coding trees with partial orders ensuring that new levels  obtained by the search engine are
 capable of extending a given fixed finite coding tree to a sub-coding tree representing a copy of $\bK$,
one is able to prove  Halpern--\Lauchli-style theorems for coding trees.
These form the pigeonhole principles of  various Milliken-style theorems for coding trees.

As the results and main ideas of the methods in  \cite{DobrinenJML20},
\cite{DobrinenH_k19}, and  \cite{Zucker20}
have been discussed in the previous section,
 we now present  the characterization of  big Ramsey degrees
in \cite{Balko7_binaryFAP}.

\begin{theorem}[\cite{Balko7_binaryFAP}]\label{thm.Balko7}
Let $\bK$ be a homogeneous structure
with finitely many relations of arity at most two
such that $\Age(\bK)=\mathrm{Forb}(\mathcal{F})$ for some finite set  $\mathcal{F}$ of  finite irreducible structures.
Then $\bK$  admits a big Ramsey structure.
\end{theorem}

Given a \Fraisse\ class
$\mathcal{K}=\mathrm{Forb}(\mathcal{F})$ with
relations of arity at most two, where $\mathcal{F}$ is a finite set of finite irreducible structures,
let $\bK$ denote an enumerated \Fraisse\ limit of $\mathcal{K}$.
Coding trees for  $\bK$ appearing in  various papers
  are all
essentially
 coding trees of $1$-types.
The proof of Theorem \ref{thm.Balko7}
uses the upper bounds of Zucker in \cite{Zucker20} as the starting point.
It then proceeds by  constructing a diagonal antichain of coding nodes which represent the structure $\bK$, with additional requirements  if there are any forbidden irreducible substructures of size three or more.
While the  exact characterization in its full generality is not
short to state,
the simpler version for the structures $\mathbf{G}_k$ include the following:
All coding nodes   $c_n\in \bfA$
code an edge with $v_m$ for some $m<n$ and
have the following property:
If  $\bfB$ is any finite graph which has the same relations over $\mathbf{G}_k\re |c_n|$ as $c_n$ does,
 then $\bfB$ has no edges.
Furthermore,
changes in the sets of  structures which are allowed to extend  a given truncation of $\bfA$ (as a level set in the coding tree)
happen as gradually as possible.
From the characterization  in
\cite{Balko7_binaryFAP},
one can make an algorithm to compute the big Ramsey degrees.

As a concrete example, we  present the exact characterization for triangle-free graphs.
In Figure 2., on the left is the beginning of   $\mathbf{G}_3$ with some fixed enumeration of the vertices as
$\{v_n:n<\om\}$.
 The $n$-th {\em coding node} in the tree
$\mathbb{S}=\mathbb{S}(\mathbf{G}_3)\sse 2^{<\om}$
 represents the $n$-th vertex $v_n$ in
$\mathbf{G}_3$, where
passing number $0$ represents a non-edge and passing number $1$ represents an edge.
 Equivalently,
 $\mathbb{S}$ is  the coding  tree of $1$-types for $\mathbf{G}_3$, as the left branch at the level of $c_n$  represents  the literal $( x \not  \hskip-.05in E v_n)$ and the right branch represents $(x E v_n)$.

%%%%%%%%%%%%%%%%%%%%%%%%%%%%
%%%%%%%%%%%%%%%%%%%%%%%%%%%
%%%%%%%%%%%%%%%%%%%%%%%%

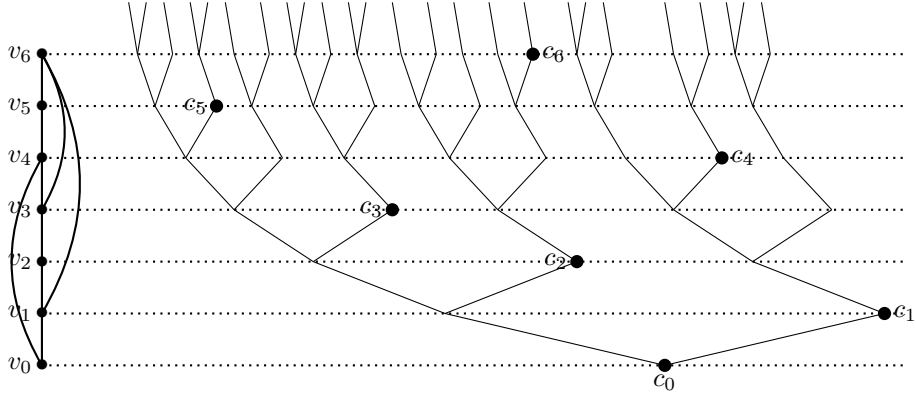
\begin{figure}[t]\label{fig.bS}
\begin{tikzpicture}[grow'=up,scale=.46]

\tikzstyle{level 1}=[sibling distance=5in]
\tikzstyle{level 2}=[sibling distance=3in]
\tikzstyle{level 3}=[sibling distance=1.8in]
\tikzstyle{level 4}=[sibling distance=1.1in]
\tikzstyle{level 5}=[sibling distance=.7in]
\tikzstyle{level 6}=[sibling distance=0.4in]
\tikzstyle{level 7}=[sibling distance=0.2in]

\node {} coordinate(t)
child{coordinate (t0)
			child{coordinate (t00)
child{coordinate (t000)
child {coordinate(t0000)
child{coordinate(t00000)
child{coordinate(t000000)
child{coordinate(t0000000)}
child{coordinate(t0000001)}
}
child{coordinate(t000001)
child{coordinate(t0000010)}
child{ edge from parent[draw=none]  coordinate(t0000011)}
}
}
child{coordinate(t00001)
child{coordinate(t000010)
child{coordinate(t0000100)}
child{coordinate(t0000101)}
}
child{ edge from parent[draw=none]  coordinate(t000011)
}
}}
child {coordinate(t0001)
child {coordinate(t00010)
child {coordinate(t000100)
child {coordinate(t0001000)}
child { edge from parent[draw=none] coordinate(t0001001)}
}
child {coordinate(t000101)
child {coordinate(t0001010)}
child { edge from parent[draw=none]  coordinate(t0001011)}
}
}
child{coordinate(t00011) edge from parent[draw=none] }}}
child{ coordinate(t001)
child{ coordinate(t0010)
child{ coordinate(t00100)
child{ coordinate(t001000)
child{ coordinate(t0010000)}
child{ coordinate(t0010001)}
}
child{ coordinate(t001001)
child{ coordinate(t0010010)}
child{ edge from parent[draw=none] coordinate(t0010011)}
}
}
child{ coordinate(t00101)
child{ coordinate(t001010)
child{ coordinate(t0010100)}
child{ coordinate(t0010101)}
}
child{   edge from parent[draw=none]coordinate(t001011)
}
}}
child{  edge from parent[draw=none] coordinate(t0011)}}}
			child{ coordinate(t01)
child{ coordinate(t010)
child{ coordinate(t0100)
child{ coordinate(t01000)
child{ coordinate(t010000)
child{ coordinate(t0100000)}
child{ edge from parent[draw=none]  coordinate(t0100001)}
}
child{ coordinate(t010001)
child{ coordinate(t0100010)}
child{edge from parent[draw=none]  coordinate(t0100011)}
}
}
child{ coordinate(t01001)
child{ coordinate(t010010)
child{ coordinate(t0100100)}
child{ edge from parent[draw=none]  coordinate(t0100101)}
}
child{edge from parent[draw=none]  coordinate(t010011)}
}}
child{ coordinate(t0101)
child{ coordinate(t01010)
child{ coordinate(t010100)
child{ coordinate(t0101000)}
child{edge from parent[draw=none]  coordinate(t0101001)}
}
child{ coordinate(t010101)
child{ coordinate(t0101010)}
child{edge from parent[draw=none]  coordinate(t0101011)}
}
}
child{  edge from parent[draw=none]  coordinate(t01011)
}}}
child{ edge from parent[draw=none]  coordinate(t011)}}}
		child{ coordinate(t1)
			child{ coordinate(t10)
child{ coordinate(t100)
child{ coordinate(t1000)
child{ coordinate(t10000)
child{ coordinate(t100000)
child{ coordinate(t1000000)}
child{ coordinate(t1000001)}
}
child{ coordinate(t100001)
child{ coordinate(t1000010)}
child{ edge from parent[draw=none] coordinate(t1000011)}
}
}
child{ edge from parent[draw=none] coordinate(t10001)
}}
child{ coordinate(t1001)
child{ coordinate(t10010)
child{ coordinate(t100100)
child{ coordinate(t1001000)}
child{edge from parent[draw=none]  coordinate(t1001001)}
}
child{ coordinate(t100101)
child{ coordinate(t1001010)}
child{edge from parent[draw=none]   coordinate(t1001011)}
}
}
child{  edge from parent[draw=none] coordinate(t10011)
}}}
child{ coordinate(t101)
child{ coordinate(t1010)
child{ coordinate(t10100)
child{ coordinate(t101000)
child{ coordinate(t1010000)}
child{ coordinate(t1010001)}
}
child{ coordinate(t101001)
child{ coordinate(t1010010)}
child{   edge from parent[draw=none]   coordinate(t1010011)}
}
}
child{   edge from parent[draw=none]  coordinate(t10101)
}}
child{ edge from parent[draw=none] coordinate(t1011)}}}
child{  edge from parent[draw=none] coordinate(t11)} };

\node[below] at (t) {$c_0$};
\node[right] at (t1) {$c_1$};
\node[left] at (t01) {$c_2$};
\node[left] at (t001) {$c_3$};
\node[right] at (t1001) {$c_4$};
\node[left] at (t00001) {$c_5$};
\node[right] at (t010101) {$c_6$};

\node[circle, fill=black,inner sep=0pt, minimum size=5pt] at (t) {};
\node[circle, fill=black,inner sep=0pt, minimum size=5pt] at (t1) {};
\node[circle, fill=black,inner sep=0pt, minimum size=5pt] at (t01) {};
\node[circle, fill=black,inner sep=0pt, minimum size=5pt] at (t001) {};
\node[circle, fill=black,inner sep=0pt, minimum size=5pt] at (t1001) {};
\node[circle, fill=black,inner sep=0pt, minimum size=5pt] at (t00001) {};
\node[circle, fill=black,inner sep=0pt, minimum size=5pt] at (t010101) {};

\draw[thick, dotted] let \p1=(t) in (-18,\y1) node (v00) {$\bullet$} -- (7,\y1);
\draw[thick, dotted] let \p1=(t1) in (-18,\y1) node (v0) {$\bullet$} -- (7,\y1);
\draw[thick, dotted] let \p1=(t01) in (-18,\y1) node (v1) {$\bullet$} -- (7,\y1);
\draw[thick, dotted] let \p1=(t001) in (-18,\y1) node (v2) {$\bullet$} -- (7,\y1);
\draw[thick, dotted] let \p1=(t1001) in (-18,\y1) node (v3) {$\bullet$} -- (7,\y1);
\draw[thick, dotted] let \p1=(t00001) in (-18,\y1) node (v4) {$\bullet$} -- (7,\y1);
\draw[thick, dotted] let \p1=(t010101) in (-18,\y1) node (v5) {$\bullet$} -- (7,\y1);

\node[left] at (v00) {$v_0$};
\node[left] at (v0) {$v_1$};
\node[left] at (v1) {$v_2$};
\node[left] at (v2) {$v_3$};
\node[left] at (v3) {$v_4$};
\node[left] at (v4) {$v_5$};
\node[left] at (v5) {$v_6$};

\draw[thick] (v0.center) to (v1.center) to (v2.center) to (v3.center);
\draw[thick] (v3.center) to (v4.center) to (v5.center);
\draw[thick] (v0.center) to [bend right] (v5.center);
\draw[thick] (v5.center) to [bend left] (v2.center);
\draw[thick] (v00.center) to [bend left] (v3.center);
\draw[thick] (v00.center) to (v0.center);
\end{tikzpicture}
\caption{Coding tree $\bS(\mathbf{G}_3)$
and the triangle-free graph represented by
its coding nodes.}
\end{figure}

%%%%%%%%%%%%%%%%%%%%%%%%%%%

Given
an antichain  $\bfA\sse \bK$,
we say that $\bfA$
is a {\em diagonal substructure} if,
letting $I$ be the set of indices of vertices in $\bfA$,
 the following hold:
(a)
For each $i\in I$,
 $v_i$ has
an edge with $v_m$ for some $m<i$;
let $m_i$ denote the least such $m$.
(b)
If $i<j$ are in $I$ with  $v_i\not \hskip-.05in E v_j$
 and  $m_j<i$,
then there is some $n\in i$ such that $v_i E v_n$ and $v_j E v_n$, and the least such $n$, denoted $n(i,j)$ is not in $I$.
(c)
For each $i, j, k, \ell \in I$ (not necessarily distinct) with $i< j$, $k< \ell$, $(i, j)\neq (k, \ell)$,  $n_j < i$, and $n_\ell< k$, we have $n(i, j)\neq n(k, \ell)$.
%For each  $J\sse I$ of size three or more such that for some $n<\min (J)$ each  vertex $v_i$, $i\in J$,  has an edge with $v_n$, letting $n_J$ be the least such $n$, for  each pairset $\{i,j\}\sse J$, there is a $k(i,j)<n_J$ such that  for all  $p\in J$, $v_p E v_{k(i,j)}$ if and only if $p\in \{i,j\}$.
Given a finite triangle-free graph $\bfA$,
the big Ramsey degree $T(\bfA)$ in $\mathbf{G}_3$ is the number of  different
 diagonal substructures representing a copy of $\bfA$.

%%%%%%%%%%%%%%%%

We conclude this section by mentioning
the exact big Ramsey degrees in the generic partial order in \cite{Balko7_PO}.
This result begins with the  upper bounds proved by \Hubicka\ in \cite{Hubicka_CS20}
and then
proceeds by
taking a diagonal antichain $D$  representing  the  generic partial order with
additional structure  of interesting levels built into $D$.
A level $\ell$ of  $D$ is {\em interesting}
if there are exactly two nodes, say $s,t$, in that level
so that
$(*)$
for exactly one  relation $\rho\in \{<,>,\bot\}$,
given any $s',t'\in D$ extending $s,t$, respectively,
$s'\rho\, t'$,
while
there is no such relation for the pair $s\re(\ell-1),t\re (\ell-1)$.
 Since an interesting level for a pair of nodes $s,t$ predetermines the relations  between any pair $s',t'$ extending $s,t$, respectively,
 passing numbers are unnecessary to the  characterization.
The big Ramsey degree of a given finite partial order  $\mathbf{P}$ is then the number of different diagonal antichains  $A\sse D$  representing $\mathbf{P}$ along with
 the order in which the interesting levels are interspersed   between the splitting levels and the  nodes in $A$.

%%%%%%%%%%%%%%%%%%%%%%%%%%%%%%%%
%%%%%%%%%%%%%%%%%%%%%%%%%%%%%%%%
%%%%%%%%%%%%%%%%%%%%%%%%%%%%%%%%

\section{Open problems
and related directions}

Section \ref{sec.theproblems} laid out the guiding questions for big Ramsey degrees.
Here we discuss some of the major open problems in big Ramsey degrees and ongoing research in  cognate areas.

\begin{problem}\label{prob.6.1}
For which SAP \Fraisse\ classes does the \Fraisse\ limit  have finite big Ramsey degrees?
\end{problem}

Subquestions are the following:
Given an SAP \Fraisse\ class with finitely many relations and a
finite set of forbidden substructures,
does its  \Fraisse\ limit have finite big Ramsey degrees?
Results in \cite{HN19}
give  evidence  for a positive answer
to this question.
For such classes with relations of arity at most two, do big Ramsey degrees always exist?
We would like a general condition on SAP classes characterizing those
 with finite big Ramsey degrees.
 We point out that
Problem \ref{prob.6.1}  in its full generality is still open for small Ramsey degrees

\begin{problem}
For results whose proofs use the method of forcing, find new proofs which are purely combinatorial.
\end{problem}

This has been done for the triangle-free graph by \Hubicka\ in \cite{Hubicka_CS20}, but new methods will be needed for $k$-clique-free homogeneous graphs for $k\ge 4$ and other such FAP classes.

The next problem has to do with topological dynamics of automorphism groups of homogeneous structures.
The work of Zucker in  \cite{Zucker19} has established a connection but not a complete correspondence yet.

\begin{problem}\label{prob.6.3}
Does every homogeneous structure with finite big Ramsey degrees also carry a big Ramsey structure?
Is there an exact correspondence, in the vein of the KPT-correspondence, between big Ramsey structures and topological dynamics?
\end{problem}

The hope in Problem \ref{prob.6.3} is to obtain as complete a dynamical understanding of big Ramsey degrees as we have for small Ramsey degrees, where a result of \cite{Zucker16} shows that given a \Fraisse\ class $\mathcal{K}$ with \Fraisse\ limit $\mathbf{K}$, then $\mathcal{K}$ has finite small Ramsey degrees if and only if the universal minimal flow of Aut$(\mathbf{K})$ is metrizable.

Finally, we mention several areas of ongoing study related to the main focus of  this paper.
Computability-theoretic
and reverse mathematical
aspects
  have been investigated
by Angl\`{e}s d'Auriac,  Cholak, Dzhafarov, Monin, and Patey.
In their treatise \cite{Cholak_etal}, they
show that the Halpern--\Lauchli\ Theorem is computably true and find  reverse-mathematical strengths for
various instances of  the product
 Milliken Theorem
and  the big Ramsey structures of the rationals and the Rado graph.
As these structures both have simply characterized big Ramsey degrees,
it will be interesting to see if different reverse mathematical  strengths emerge
 for
structures
such as  the triangle-free homogeneous graph or the generic partial order.

Extending Harrington's forcing proof to
the uncountable realm,
 Shelah in
 \cite{Shelah91}  showed that  it is consistent, assuming certain large cardinals, that  the Halpern--\Lauchli\ Theorem  holds
for trees $2^{<\kappa}$, where $\kappa$ is
a measurable cardinal.
D\v{z}amonja, Larson, and Mitchell applied a slight modification of his theorem to characterize the   big Ramsey degrees for the $\kappa$-rationals and the $\kappa$-Rado graph.
Their characterizations have as their basis the characterizations of
 Devlin and Laflamme--Sauer--Vuksanovic for the rationals and Rado graph, respectively,
but also involve  well-orderings of each level of the tree $2^{<\kappa}$, necessitated by  $\kappa$ being uncountable.
The field of big Ramsey degrees for uncountable homogeneous structures  is still  quite open, but
the fleshing out of the Ramsey theorems on trees  of uncountable height has seen
some recent work   in \cite{DH16}, \cite{DH18}, and \cite{Zhang17}.

The next  problem
comes  from a general question in \cite{KPT05}.

\begin{problem}
Develop infinite dimensional Ramsey theory on spaces of copies of a  homogeneous structure.
\end{problem}

For a set $N\sse\om$, let $[N]^{\om}$ denote the set of all infinite subsets of $N$, and note that $[\om]^{\om}$ represents the Baire space.
The infinite dimensional
Ramsey theorem
of  Galvin and Prikry (\cite{Galvin/Prikry73})
says that
 given any Borel subset $\mathcal{X}$ of the Baire space,
there is an infinite set $N$
such that  $[N]^{\om}$ is either contained in $\mathcal{X}$ or is disjoint from $\mathcal{X}$.
Ellentuck's theorem in  \cite{Ellentuck74}
found optimality in terms of sets with the property of Baire with respect to
a finer topology.
The question in \cite{KPT05} asks for extensions of these theorems to
subspaces of $[\om]^{\om}$, where each infinite set represents
 a copy of some fixed homogeneous structure.
A Galvin--Prikry style theorem
 for spaces of copies of the  Rado graph
has been proved by the author in
\cite{DobrinenRado19}.
By a comment of Todorcevic in Luminy in 2019, the infinite dimensional Ramsey theorem should ideally also recover  exact big Ramsey degrees.
Such a theorem is being written down by the author for structures satisfying SDAP$^+$ with relations of arity at most two.
This is one instance where coding trees are necessitated to be diagonal  in order for the infinite dimensional Ramsey theorem to directly recover exact big Ramsey degrees.

We close by mentioning that  structural Ramsey theory has  been central in investigations of ultrafilters which are  relaxings  of Ramsey ultrafilters in the same way that big Ramsey degrees are relaxings  of  Ramsey's theorem.
An exposition of recent work appearing in
\cite{DobrinenSEALS}
will give the reader  yet another view of the power of  Ramsey theory.
\vskip.1in

%------
% Insert acknowledgments and information
% regarding funding at the end of the last
% section, i.e., right before the bibliography.
%------

\noindent\bf Acknowledgements. \rm  The author is grateful to  her PhD advisor Karel Prikry, and  to
 Richard Laver,  Norbert Sauer,  and Stevo Todorcevic who have mentored  her in areas related to this paper.
She thanks Jan \Hubicka,
Dragan Ma\v{s}ulovi\'{c},
 Mat\v{e}j Kone\v{c}n\'{y},
Lionel Nguyen Van Th\'{e},  and Andy Zucker for valuable  feedback  on drafts of this paper.

%------
% Insert the bibliography.
%------

%\bibliographystyle{amsplain}

\end{document}